\documentclass[10pt]{amsart}
\usepackage{amsmath,amssymb,amsthm,mathtools}
\usepackage{fullpage, xcolor}
\usepackage{hyperref}
\usepackage{booktabs}

\newtheorem{thm}{Theorem}[section]
\newtheorem{lem}[thm]{Lemma}
\newtheorem{prop}[thm]{Proposition}
\newtheorem{cor}[thm]{Corollary}
\theoremstyle{definition}
\newtheorem{defs}[thm]{Definition}
\theoremstyle{remark}
\newtheorem{rem}[thm]{Remark}
\numberwithin{equation}{section}

\newcommand{\sine}{\mathsf{sine}}

\newcommand{\Var}{\mathrm{Var}}

\newcommand{\R}{{\mathbb R}}
\newcommand{\Z}{{\mathbb Z}}
\newcommand{\mbbC}{\mathbb C}
\newcommand{\mbbE}{\mathbb E}
\newcommand{\mbbP}{\mathbb P}
\newcommand{\mcN}{\mathcal N}
\newcommand{\mcC}{\mathcal C}
\newcommand{\mcH}{\mathcal H}

\newcommand{\mcS}{\mathcal S}
\newcommand{\mcX}{\mathcal X}
\newcommand{\mcJ}{\mathcal J}
\newcommand{\mcR}{\mathcal R}
\newcommand{\mbfK}{\mathbf K}
\newcommand{\1}{\mathbf 1}
\newcommand{\Cov}{\mathrm{Cov}}
\DeclareMathOperator{\Tr}{Tr}
\newcommand{\supp}{\mathrm{supp}}
\newcommand{\law}{\mathrm{law}}
\newcommand{\tr}{\mathrm{tr}}

\newcommand{\mrmconst}{\mathrm{const}}
\newcommand{\defi}{\coloneqq}
\newcommand{\diff}{\mathop{}\!\mathrm{d}}
\newcommand{\comma}{\,\,\mathrm{,}\;\,}

\newcommand{\fstop}{\,\,\mathrm{.}}

\title{\bf The singular Gaussian scaling limit for a tagged particle in the infinite-dimensional Dyson model at $\beta=2$}
\author{Hirofumi Osada}
\address{Chubu University
1200 Matsumoto-cho, Kasugai-shi
Aichi 487-8501, Japan}
\email{osada@fsc.chubu.ac.jp}

\author{Kohei Suzuki}
\address{Department of Mathematical Sciences, Durham University, South Road, Durham
DH1 3LE, United Kingdom}
\email{kohei.suzuki@durham.ac.uk}

\thanks{Hirofumi Osada: Chubu University
1200 Matsumoto-cho, Kasugai-shi
Aichi 487-8501, Japan
\texttt{osada@fsc.chubu.ac.jp}. \hfill
\\
\indent Kohei Suzuki: Durham University, Durham, United Kingdom;
\texttt{kohei.suzuki@durham.ac.uk}}

\subjclass[2020]{Primary 60K35; Secondary 60B20, 82C22}
\keywords{Dyson Brownian motion, sine process, tagged particle, central limit theorem, determinantal point process}

\newcommand{\mcE}{\mathcal E}

\begin{document}

\begin{abstract}
\fontsize{10pt}{12pt}\selectfont
 We prove a logarithmic central limit theorem for a tagged particle in infinite-
dimensional Dyson Brownian motion at inverse temperature $\beta=2$, starting from the
sine$_2$ point process. The fluctuation is of order $\sqrt{\log t}$, which is subdiffusive and substantially slower than the polynomial fluctuations of classical short-range interacting particle systems,
reflecting the rigidity induced by long-range logarithmic repulsion. We further identify the scaling limit as the singular Gaussian field corresponding to the zero-Hurst parameter, which provides a realisation of the singular Gaussian field as a scaling limit of a tagged particle raised by Dürr, Goldstein and Lebowitz (Comm. Pure Appl. Math. 38 (1985)).

\end{abstract}

\maketitle

% \purple{\bf Cite two papers by Dandekar Krapivsky Mallick} 
% \\
% \purple{\bf "https://arxiv.org/pdf/2409.06881" and https://arxiv.org/abs/2212.05583}

\section{Introduction}

A basic question in the theory of interacting particle systems is
how interactions affect the long-time behaviour of a tagged particle.
In one-dimensional systems of infinitely many diffusion processes where particles cannot pass one
another, the motion of a tagged particle is constrained by
its neighbours. This restriction can lead to fluctuations
on a scale smaller than the usual diffusive scale $t^{1/2}$.
%even when the particles almost surely never collide.
For one-dimensional continuum systems, Harris~\cite{Harris1965} studied infinitely many Brownian motions with Poisson initial positions, where particles move independently until they collide and they instantaneously reflect each other after collision. We call such an interaction  {\it collision-reflection interaction}.  
Under this setting, he proved the central limit theorem (CLT) with the $t^{1/4}$ displacement scale for a tagged particle  converging to the standard Gaussian random variable.

Beyond the CLT of one-time marginal, one may study the multi-time scaling limit of a tagged particle. 
For independent exponential initial
particles and independent, identically distributed centred velocities with finite, nonzero first absolute
moment, Spitzer~\cite{Spitzer1969} proved a Brownian scaling limit for the tagged particle in an infinite system of collision–reflection interaction with uniform free motion of constant velocity. 
The celebrated paper by D\"urr--Goldstein--Lebowitz~\cite[Theorem 7.1 and Example 8.4]{DurrGoldsteinLebowitz1985} considerably strengthened this result to the case where driving noises are including fractional Brownian motions with all Hurst parameter~$0<H  \le 1$ (the case $H=1$ is the uniform free motion with random velocity in~Spitzer~\cite{Spitzer1969}). As the scaling limit with fluctuation order $t^{H/2}$,  they obtained  the fractional Brownian motion with Hurst parameter $H/2$ for every $0<H \le 1$. 
In~\cite[Section~4, final remark following
(4.12)]{DurrGoldsteinLebowitz1985}, 
they raised the question whether the
singular Gaussian field corresponding to the Hurst parameter~$H=0$ could arise as a scaling limit of
their process. 

\smallskip
In this paper, we show that this singular Gaussian field arises as 
the multi-time scaling limit of a tagged particle in the
stationary infinite-dimensional Dyson Brownian motion at the inverse temperature~$\beta=2$, where the interaction is {\it long-range (logarithmic) and repulsive}.
Due to the long-range repulsive interaction, large scale density fluctuations are strongly suppressed,   and the variance of the tagged-particle displacement grows only logarithmically in time, which is subdiffusive, and remarkably slower than the aforementioned local interaction cases. This logarithmic variance growth has been originally discovered by Spohn~\cite{Spo87}, while the corresponding CLT and the multi-time scaling limit are left open. 
As a main result, we prove that the multi-time scaling limit of a tagged particle with the logarithmic scale~$\sqrt{\log t}$ exists as the convergence in finite-dimensional distribution. Furthermore,  we  identify the limit with the singular Gaussian with~$H=0$.

\bigskip
\paragraph{\bf Main results} The infinite-dimensional Dyson Brownian motion (DBM) is the bulk scaling limit of finite DBM and has the
sine point process $\sine_2$ as its invariant distribution. Its determinantal
and Markovian structure was developed in
\cite{KatoriTanemura2010,KatoriTanemura2011,OsadaTanemura2016}. The realisation of DBM as a weak solution of the associated infinite-dimensional SDE
 was proven in~\cite{Osada2012}; strong existence and pathwise uniqueness
in suitable initial conditions were established in~\cite{Tsai2016, OsadaTanemura2020}.

To be more precise, we consider an infinite system of Brownian particles on $\R$
with logarithmic repulsion at inverse temperature $\beta=2$.
Let $(X_i(t))_{i\in\Z,t\ge0}$ be continuous coordinate processes
satisfying
$$
    \cdots<X_{-1}(t)<X_0(t)<X_1(t)<\cdots
    \qquad\text{for all }t\ge0
$$
almost surely, and solving the infinite-dimensional stochastic
differential equations
\begin{align} \label{e:DSDE}
    \diff X_i(t)
    &=
    \lim_{r\to\infty}
    \sum_{\substack{j\ne i\\ |X_i(t)-X_j(t)|<r}}
    \frac{\diff t}{X_i(t)-X_j(t)}
    +\diff B_i(t)\comma
    \qquad i\in\Z\fstop
\end{align}
Here $(B_i)_{i\in\Z}$ are independent standard one-dimensional
$(\mathcal F_t)$-Brownian motions defined on a filtered probability
space $(\Omega,\mathcal F,(\mathcal F_t)_{t\ge0},\mbbP)$ satisfying
the usual conditions. We write $\mbbE$ for expectation with respect
to~$\mbbP$.
The notation $Z\sim\mu$ means that the random element $Z$ has
law $\mu$ under $\mbbP$. For deterministic quantities, $a\sim b$
denotes the asymptotic equivalence, i.e.,  $a/b\to1$ in the specified
limit.

The point configuration associated to $(X_i(t))_{i \in \Z, t \ge 0}$ is
$$
    \xi_t\defi\sum_{i\in\Z}\delta_{X_i(t)}\fstop
$$
We take the stationary solution, i.e., $(\xi_t)_{t\ge0}$
is stationary such that
$$
    \xi_t\sim\sine_2
    \qquad\text{for every }t\ge0\fstop
$$
Here $\sine_2$ is the law of the determinantal point process
on $\R$, with respect to Lebesgue measure, whose correlation
kernel is
$$
    K_{\sin}(x-y)\defi
    \frac{\sin\pi(x-y)}{\pi(x-y)}\comma
    \qquad K_{\sin}(0)\defi1\fstop
$$
We call $(\xi_t)_{t\ge0}$ the stationary unlabelled Dyson
Brownian motion and the coordinates $(X_i(t))_{i\in\Z}$ keeps
the particle labels.

The drift of~the SDE~\eqref{e:DSDE} is the repulsive force associated with the logarithmic
pair potential $-\log|x|$.
This
long range repulsive force  makes the equilibrium configuration rigid, and  it produces the logarithmic fluctuation of a tagged-particle  studied below. We choose the labelling so that
$X_0(0)\le0<X_1(0)$ and study the tagged particle $X_0(t)$.

Our main results establish the logarithmic CLT and the corresponding multi-time scaling limit. 
\begin{thm}[CLT for the tagged particle, Corollary~\ref{cor:tracer-CLT}]\label{t:tracer-CLT}
As $t\to\infty$, the following convergence in distribution holds:
$$
\frac{\pi\,(X_0(t)-X_0(0))}{\sqrt{\log t}}
\ \Longrightarrow\ \mcN(0,1)\fstop
$$
 
\end{thm}

\begin{thm}[Multi-time scaling limit, Theorem~\ref{thm:tagged-process-clt}]\label{thmI:tagged-process-clt}
Let
$$
\Delta(t)\defi X_0(t)-X_0(0)\comma \qquad
\mcX_\alpha(u)\defi \frac{\pi\Delta(u/\alpha)}{\sqrt{\log(1/\alpha)}}\fstop
$$
Then, for every finite collection $u_1,\ldots,u_n>0$,
$$
(\mcX_\alpha(u_1),\ldots,\mcX_\alpha(u_n))
\Longrightarrow
(G(u_1),\ldots,G(u_n))\comma \qquad \alpha\downarrow0\comma
$$
where $G$ is the centred Gaussian field with covariance
\begin{equation}\label{eqI:homogenized-current-covariance}
\mbbE[G(u)G(v)]
=\frac12+\frac12\1_{\{u=v\}}
=\begin{cases}
1,&u=v,\\[1mm]
\dfrac12,&u\ne v.
\end{cases}
\end{equation}
\end{thm}
This covariance formula appears in the limit $H \to 0$ of the Gaussian process $G_H$ with Hurst parameter~$H$:
$$\mbbE[G_H(u)G_H(v)]
=\frac12(u^{2H} + v^{2H} - |u-v|^{2H}) \xrightarrow{H \to 0} \frac12 + \frac12\1_{\{u=v\}}\fstop$$

\smallskip
\paragraph{\bf Ideas of the proofs.} 
A key observation regarding the logarithmic scale has been already made by 
Spohn~\cite{Spo87}:  for the equilibrium infinite Dyson model at
$\beta=2$, the mean-square displacement of a tagged particle grows only
logarithmically:
\begin{align} \label{e:SP}
    \mbbE\bigl[|X_0(t)-X_0(0)|^2\bigr]
    \sim \frac{1}{\pi^2}\log t\comma \qquad t\to\infty\fstop
\end{align}
The proof in~\cite{Spo87} is based on the estimate of the current $J(t)$ defined by ``number of particles that cross $0$ from left to right" minus ``number of particles that cross $0$ from right to left" during the time interval $[0,t]$ in~\cite[pp.~671--673]{Spo87}.
This description expresses the current $J(t)$ as the net transport across the origin.
However, recalling the fact that each one-dimensional Brownian motion can recross a level infinitely many time in any compact time interval, if individual crossing events are counted literally, this  definition ends up with the indefinite expression $``\infty-\infty"$. 
In this paper, we define the net current in the
infinite system as the  difference of finite index
\begin{equation}\label{d:RDJ}
 J(t)\defi I_0(0)-I_0(t)\comma
 \qquad I_u(t)\defi\sup\{i\in\Z:X_i(t)\le u\}\fstop
\end{equation}
See Definition~\ref{d:Cur}. 
Lemma~\ref{lem:finite-sum} identifies it with the sum of the
signed changes of side between times $0$ and $t$, with only
finitely many nonzero particle contributions.

To justify the variance calculation, we approximate this
pathwise current by compactly supported linear statistics:
$$
 f_R(x)\defi\chi(x/R)\1_{\{x>0\}}\comma
 \qquad
 Z_R(t)=\int_\R f_R(x)
       \bigl(\xi_t(\diff x)-\xi_0(\diff x)\bigr)\fstop
$$
Here $\chi$ is a fixed compactly supported smooth cut-off with
$0\le\chi\le1$, whch is equal to one near the origin. 
The space-time determinantal structure of the stationary DBM $(\xi_t)_{t \ge 0}$ established in
\cite{KatoriTanemura2010,KatoriTanemura2011,OsadaTanemura2016}
gives the covariance formulas in Proposition~\ref{prop:cov-FT} and Lemma~\ref{l:WCV}, which
allow us to control the distant boundary effect.
Proposition~\ref{prop:direct-current-L2} proves that $J(t)$ is
centred and square-integrable, with the error estimate
$$
 \mbbE|Z_R(t)-J(t)|^2\le C\frac{t}{R}\comma
 \qquad t>0,\quad R\ge1\fstop
$$
This $L^2$ convergence justifies passage from the variance of the compact-statistic $Z_R(t)$
to the one of the current~$J(t)$ in
Proposition~\ref{prop:Var-J}. Theorem~\ref{thm:main} then gives
$$
 \Var(J(t))=\frac{1}{\pi^2}\log t+O(1)\fstop
$$
These steps supply the limiting argument underlying Spohn's
test-function approximation in~\cite[pp.~672--673]{Spo87};
see Remark~\ref{r:CS} for the detailed comparison.

In section~\ref{sec:CLT}, we prove the CLT. We take $X_t\defi Z_{t^2}(t)$ and prove the approximation bound~
$\|X_t-J(t)\|_{L^2}=O(t^{-1/2})$, and then
Lemma~\ref{lem:CLT-reduction} reduces the problem to proving
the CLT for $X_t$ in place of $J(t)$.
Proposition~\ref{prop:fredholm} expresses the moment generating
function of $Z_R(t)$ as a Fredholm determinant.
Taking its logarithm and comparing coefficients gives the
cumulant formula in Proposition~\ref{prop:cumulant-trace}.
For each fixed $m\ge3$, the $m$th cumulant is a finite weighted
sum of integrals of the form
$$
 \int_{\R^\ell}
 \prod_{j=1}^{\ell}
 \Bigl[f_R(x_j)^{p_j}
 K_{\sin}(\tau_j,x_j;\tau_{j+1},x_{j+1})\Bigr]
 \diff x_1\cdots\diff x_\ell\comma
$$
where $K_{\sin}(\tau_j,x_j;\tau_{j+1},x_{j+1})$ denotes the space-time determinantal kernel of~$(\xi_{t})_{t \ge 0}$.  
Here $1\le\ell\le m$, the positive integers $p_j$ sum to $m$,
and $\tau_j\in\{0,t\}$, with
$(\tau_{\ell+1},x_{\ell+1})=(\tau_1,x_1)$.
Thus the kernel factors form a closed cycle.
The coefficients include the signs from subtracting the
time-$0$ statistic. We estimate these integrals in two groups.

First, collect all terms using only one observation time,
summing over every length $\ell$ and composition of~$m$
with the coefficients from Proposition~\ref{prop:cumulant-trace}.
Their total contribution is given by
$$
 \kappa_m(\mcN_t(f_R))
 +(-1)^m\kappa_m(\mcN_0(f_R))\fstop
$$
The cyclic integrals appearing in the expansion of one-time cumulant can  be of order \(R\).  Their signed sum is however of order \(O_m(\log R)\)  due to the cancellation, see~Proposition~\ref{prop:sine-cumulant-bound}.   Consequently, the contribution from the two time points also has the logarithmic bound above by $C_m\log R$, where $C_m$ is a constant depending only on~$m$.

The Remaining integral uses both times points. 
Lemma~\ref{lem:two-interfaces} with the power estimate
\eqref{eq:comm-power} gives a bound its absolute value by $m^2d_R^2$, viz., 
$$
 d_R^2\defi
 \iint_{\R^2}|f_R(x)-f_R(y)|^2
 \left|\frac{\sin\pi(x-y)}{\pi(x-y)}\right|^2
 \diff x\diff y\fstop
$$
Lemma~\ref{lem:comm-logR} provides~$d_R^2\le C\log R$.

For fixed $m$, the number of terms and their coefficients depend
only on $m$. Summing the preceding bounds proves the two-time
case of Proposition~\ref{prop:cumulant-bound}:
$$
 |\kappa_m(Z_R(t))|\le C_m\log R\comma
 \qquad R\ge2,\quad t>0\fstop
$$
With $R=t^2$, Theorem~\ref{thm:main} and
Lemma~\ref{lem:CLT-reduction} give
$\Var(X_t)\sim\pi^{-2}\log t$, while homogeneity of cumulants gives
$$
 \left|\kappa_m\!\left(\frac{X_t}{\sqrt{\log t}}\right)\right|
 \le C_m(\log t)^{1-m/2}\longrightarrow0\comma
 \qquad m\ge3\fstop
$$
The method of moments together with Spohn's variance formula yields the Gaussian limit, see Theorem~\ref{thm:CLT}.

In order to go from the current estimate to the tagged particle estimate,
we prove Lemma~\ref{lem:N-index}, which gives the identity
$J(t)=N_t(X_0(t))$, where $N_t(x)$ counts particles between
$0$ and $x$ with a negative sign when $x<0$.
The equilibrium density is one, thus $\mbbE[N_t(x)]=x$.
By using Lemma~\ref{lem:number-variance} and Chebyshev's inequality, we can control the term
$N_t(n)-n$ for every integer~$n$.  
We then prove in
Proposition~\ref{prop:tracer-minus-current} that 
$$
 \frac{X_0(t)-X_0(0)-J(t)}{\sqrt{\log t}}
 \xrightarrow{\mbbP}0\fstop
$$
The CLT for the current therefore implies the CLT for the tagged-particle, see
Corollary~\ref{cor:tracer-CLT} for the details.

To prove the multi-time scaling limit, we apply the same argument to 
linear combinations of currents at times $u_1/\alpha,\ldots,u_n/\alpha$,
where $u_j>0$ and $\alpha\downarrow0$.
By Lemma~\ref{lem:simultaneous-compact-approx}, we can replaces these
currents simultaneously by linear statistics~$Z_{R_\alpha}(u_j/\alpha)$ with the
common cut-off $R_\alpha=\alpha^{-2}$.
The bound for multi-time cumulants in Proposition~\ref{prop:cumulant-bound}
applies and  the higher cumulants
vanish after normalisation by $\sqrt{\log(1/\alpha)}$.

For the covariance,  by writing
$V(t)\defi\Var(J(t))$ and  using the stationarity of the current increments, we obtain 
$$
 \Cov(J(s),J(t))
 =\frac12\bigl(V(s)+V(t)-V(|t-s|)\bigr)\fstop
$$
For $s=u/\alpha$ and $t=v/\alpha$ with  $u\ne v$,
all these variances have the same leading logarithmic scale due to
Theorem~\ref{thm:main}. Thus, after normalisation by
$\pi/\sqrt{\log(1/\alpha)}$, the limiting variance is $1$
and the limiting covariance at distinct positive times is $1/2$,
as given  in Lemma~\ref{lem:current-covariance}.
By the Cram\'er--Wold theorem, we obtain the 
Gaussian scaling limit in convergence of finite-dimensional distribution, see
Theorem~\ref{thm:current-process-clt}.
Applying Proposition~\ref{prop:tracer-minus-current} at finitely many observation times transfers this limit to
the tagged-particle displacements. This proves
Theorem~\ref{thm:tagged-process-clt}.

\smallskip
\paragraph{\bf Further literature} Here we supply related literature in the case where driving noises are continuous in time. 
The infinite Atlas model considers the lowest ranked Brownian particle receives a
constant positive drift. Starting from the stationary
distribution of the gaps, Dembo--Tsai~\cite[Corollary~1.4(a)]
{DemboTsai2017} proved that the lowest ranked particle,
under the time scaling $t\mapsto Tt$ and spatial scaling
$T^{-1/4}$, converges to a constant multiple of fractional
Brownian motion with Hurst parameter $1/4$. 
Osada~\cite[Theorem 2.1]{Osada2016SelfDiffusion}
proved, in the Dirichlet form setting, that non-collision forces the
self-diffusion constant to vanish. This rules out a nondegenerate diffusive limit but does not identify an exact slower
fluctuation scale. 

Dandekar–Krapivsky–Mallick~\cite{DandekarKrapivskyMallick2023} studied Brownian particles with repulsive Riesz potential proportional to $|x-y|^{-s}$ with $s>0$. By using hydrodynamical arguments for fluctuations, they obtained that the standard-deviation of the current and the tagged-particle  have the scale $t^{s/[2(1+s)]}$ for $0<s<1$ and $t^{1/4}$ for $s>1$. For $0<s<1$, their equilibrium two-time covariance calculations match the covariance of the fractional Brownian motions, although  a corresponding functional limit theorem is left open.

For the Dyson gas, Dandekar–Krapivsky–Mallick~\cite{DandekarKrapivskyMallick2024} used a macroscopic hydrodynamic computation to rederive Spohn’s leading logarithmic variance asymptotic and derive a heuristic computation of Gaussian current fluctuations with formal higher-cumulant asymptotics. Interestingly, their analysis applies beyond $\beta=2$ while parts of cumulant computations are left heuristic, e.g., it relies on a saddle-point approximation to the current’s cumulant generating function [Eqs. (20)–(21)], without controlling the approximation error sufficiently to establish a CLT for the underlying particle dynamics. For the stationary infinite Dyson system at $\beta=2$, we use a microscopic approach from the SDE perspective and give a pathwise probabilistic derivation of the variance asymptotic and rigorously prove the current and tagged-particle CLTs, together with the multi-time Gaussian scaling limit.

\subsection*{Declaration of AI Use}
The main theorems and  the proof ideas of  the cumulants estimates by using the space-time determinantal structure of Dyson model and the Fourier transform were due to the authors. The details of computations, particularly higher-order cumulant computations and the related literature search for operator algebras were assisted by~interactive discussions with GPT~5.5,  5.6 and 6. We also used them for language polishing.   All the mathematical arguments were independently verified by the authors who take the full responsibility for the results.

\subsection*{Acknowledgement}
We are grateful to 
Paul Krapivsky and Kirone Mallick for letting us know about their papers~\cite{DandekarKrapivskyMallick2023, DandekarKrapivskyMallick2024}. 
HO was supported by JSPS KAKENHI Grant Numbers JP21H04432 and JP24KK0060. KS was supported by the Ito-Arakawa Foundation Grant from the Humboldt Association of Japan.

\subsection*{Notation}
The \emph{Schwartz space} $\mcS(\R)$ consists of all smooth
functions $f:\R\to\R$ such that, for every pair of non-negative
integers $m,n$,
$$
 \sup_{x\in\R}(1+|x|)^m\,|f^{(n)}(x)|<\infty\fstop
$$
Thus $f$ and all its derivatives decay faster than any inverse
power of $|x|$ as $|x|\to\infty$.
We use the following Fourier transform and inversion conventions:
\begin{equation}\label{eq:FT}
\widehat f(k)=\int_{\R} e^{-ikx} f(x)\diff x\comma \qquad
f(x)=\frac{1}{2\pi}\int_{\R} e^{ikx}\widehat f(k)\diff k\fstop
\end{equation}
For functions $f,g$ on $\R$, the symbol $*$ denotes convolution:
$$
 (f*g)(x)\defi\int_\R f(x-y)g(y)\diff y\comma
$$
whenever the integral is absolutely convergent. With the convention
in~\eqref{eq:FT}, we have
$$
 \widehat{f*g}(k)=\widehat f(k)\widehat g(k)\comma
 \qquad f,g\in L^1(\R)\fstop
$$
An asterisk in a superscript, as in $T^*$, denotes the adjoint
of an operator.

For real-valued Schwartz functions $f$ and $g$, we use the bilinear form of Plancherel's identity:
$$
\int_{\R} f(x)g(x)\diff x
=\frac{1}{2\pi}\int_{\R}\widehat f(k)\widehat g(-k)\diff k\fstop
$$
%A real-valued function need not have a real-valued Fourier transform. Nevertheless,
For real $f$, we have $\widehat f(-k)\widehat f(k)=|\widehat f(k)|^2$, so
$$
\int_{\R}f(x)^2\diff x
=\frac{1}{2\pi}\int_{\R}|\widehat f(k)|^2\diff k\fstop
$$
%Thus no explicit complex conjugates are needed in our covariance formulae, but
%squared moduli remain essential in the variance formulae.
At density $1$, the sine kernel is
\begin{equation}\label{eq:sine-kernel}
K_{\sin}(x-y)=\frac{\sin\pi(x-y)}{\pi(x-y)}=\frac{1}{2\pi}\int_{-\pi}^{\pi} e^{ik(x-y)}\diff k\comma
\end{equation}
i.e.\ the Fourier multiplier is the Fermi projector $\1_{\{|k|\le \pi\}}$:
$$
\widehat{Pf}(k)
=\widehat{K_{\sin}*f}(k)
=\1_{\{|k|\le\pi\}}\widehat f(k)\comma
\qquad f\in L^2(\R)\fstop
$$

We shortly write $L^p$ 
for $L^p(\Omega,\mathcal F,\mbbP)$, for $1\le p\le\infty$.
In particular,
$$
\|Y\|_{L^2}\defi\bigl(\mbbE[|Y|^2]\bigr)^{1/2}\comma
\qquad
Y_n\xrightarrow{L^2}Y
\quad\Longleftrightarrow\quad
\mbbE[|Y_n-Y|^2]\longrightarrow0\fstop
$$
Spatial $L^p$ spaces, such as $L^2(\R)$, and the $L^p$ norms of
deterministic test functions are taken with respect to Lebesgue measure,
unless another measure is specified.

\section{Logarithmic asymptotics for the variance of the current}
We discuss the logarithmic asymptotic of variance  for the current
across the origin:
\begin{align}\label{e:VF}
\Var(J(t))=\frac{1}{\pi^2}\log t+O(1)\comma
\qquad t\to\infty\fstop
\end{align}
This statement was already provided in~\cite{Spo87}.  As we need several details skipped there  for the sake of the multi-time scaling limit  in the later sections, we provide detailed discussions regarding~\eqref{e:VF} in this section. See Remark~\ref{r:CS} for further comparisons with~\cite{Spo87}.

The proof has three main steps. First, we compute the covariance of
compactly supported linear statistics and its Fourier transform.
Next, we show that cutting off the half-line far from the origin gives
an approximation of the pathwise current in the $L^2$ sense. Finally, the variance formula by Fourier transform provides a frequency
band whose contribution is logarithmic.
The same compact half-line statistics will be used in the cumulant
argument of Section~\ref{sec:CLT}.

\subsection{The extended sine kernel and the current}
\label{subsec:extended-sine-current}
Let $(\xi_t)_{t\ge0}$ be the stationary unlabelled Dyson diffusion
at inverse temperature $\beta=2$ and density $1$, starting with  invariant
distribution $\sine_2$. We first specify the correlation kernel
used throughout the paper.

\begin{defs}[Extended sine kernel]\label{d:extended-sine-kernel}
For $s,r\ge0$ and $x,y\in\R$, the {\it extended sine kernel with
density $1$} is defined as
\begin{equation}\label{eq:extended-sine-kernel}
 K_{\sin}(s,x;\,r,y)\defi
 \begin{cases}
 \displaystyle\frac1\pi\int_0^\pi
 e^{(r-s)k^2/2}\cos\!\bigl(k(x-y)\bigr)\diff k,
 &s\le r \comma\\[2mm]
 \displaystyle-\frac1\pi\int_\pi^\infty
 e^{-(s-r)k^2/2}\cos\!\bigl(k(x-y)\bigr)\diff k,
 &s>r \fstop
 \end{cases}
\end{equation}
The normalisation used here is the same one used in~\cite[Eq.~(1.4)]{KatoriTanemura2011}.
When~$s=r$ it comes down to the usual static sine kernel
in~\eqref{eq:sine-kernel} with value $1$ when $x=y$.
Equivalently, with the convention~\eqref{eq:FT},
\begin{equation}\label{eq:ext-sine-FT}
 K_{\sin}(s,x;\,r,y)
 =\frac1{2\pi}\int_{\R}e^{ik(x-y)}
       \kappa_{s,r}(k)\diff k\comma
\end{equation}
where the spatial Fourier multipliers are
\begin{equation}\label{eq:kappa-KT}
 \kappa_{s,r}(k)\defi
 \begin{cases}
 e^{(r-s)k^2/2}\1_{\{|k|\le\pi\}},&s\le r\comma \\
 -e^{-(s-r)k^2/2}\1_{\{|k|>\pi\}},&s>r \fstop
 \end{cases}
\end{equation}
\end{defs}

The stationary Dyson Brownian motion has the space-time determinantal representation
with this kernel; see~\cite[Theorem~2.2]{OsadaTanemura2016}, after
the deterministic space-time scaling to density $1$. In particular,
on any finite collection of time layers, its $n$-point correlation
density with respect to counting measure in time and Lebesgue
measure in space is
\begin{equation}\label{eq:extended-sine-correlations}
 \rho^{(n)}(z_1,\ldots,z_n)
 =\det\!\bigl[K_{\sin}(s_i,x_i;\,s_j,x_j)\bigr]_{i,j=1}^n\comma
 \qquad z_i=(s_i,x_i)\fstop
\end{equation}
The associated path law in the continuous path space over the configuration space (i.e., the space of locally finite point measures over~$\R$ endowed with the vague topology)  is invariant under spatial translations
and reflections, and is time reversible; see
\cite[Theorem~2.6]{KatoriTanemura2011}.

Furthermore, particles do not collide within finite time almost surely, see~Osada~\cite{Osada2004}, and also Assiotis--Suzuki~\cite{AssiotisSuzuki2026}.  Choose an increasing labelling of the initial configuration and follow these labels along the non-colliding paths. Thus
$$
\xi_t=\sum_{i\in\Z}\delta_{X_i(t)}\comma \qquad
\cdots<X_{-1}(t)<X_0(t)<X_1(t)<\cdots\fstop
$$
The labelled particles satisfy the infinite-dimensional stochastic differential equation
\begin{align} \notag %\label{DBM}
\diff X_i(t)
=\lim_{r\to\infty}\sum_{\substack{j:j\ne i\\ |X_i(t)-X_j(t)|<r}}
\frac{\diff t}{X_i(t)-X_j(t)}+\diff B_i(t)\comma \qquad i\in\Z\comma
\end{align}
where $(B_i)_{i\in\Z}$ are independent standard Brownian motions and $\xi_0\sim\sine_2$, see~\cite{Osada2012, Tsai2016, OsadaTanemura2020}.

For $u\in\R$ and $t\ge0$, define
\begin{align}\label{:22}&
I_u(t)\defi\sup\{i\in\Z:X_i(t)\le u\}\fstop
\end{align}
For each fixed $t\ge0$, almost surely the configuration is locally finite and has infinitely many points in both directions. Hence $X_i(t)\to\pm\infty$ as $i\to\pm\infty$, respectively, and $I_u(t)\in\Z$ is well-defined for every $u\in\R$.

\begin{defs}[Current] \label{d:Cur}
For $u\in\R$ and $t\ge0$, define the current across level $u$ by
\begin{equation}\notag %\label{eq:J-def}
J(u,t)\defi I_u(0)-I_u(t)\in\Z\fstop
\end{equation}
We write $J(t)\defi J(0,t)$ for the current across the origin. With this sign convention, the movement from left to right contributes positively to the current.
\end{defs}

\begin{lem}\label{lem:finite-sum}
For each fixed $u\in\R$ and $t\ge0$, almost surely,
$$
J(u,t)
=\lim_{L\to\infty}\sum_{i\in\Z\cap[-L,L]}
\bigl(\1_{\{X_i(t)>u\}}-\1_{\{X_i(0)>u\}}\bigr)\comma
$$
and only finitely many summands are non-zero.
\end{lem}

\begin{proof}
The ordering of the particles gives
$$
\1_{\{X_i(s)>u\}}=\1_{\{i>I_u(s)\}}\comma \qquad s\in\{0,t\}\fstop
$$
If $I_u(0)<I_u(t)$, the non-zero summands are
those with labels
$$
 i=I_u(0)+1,\ldots,I_u(t)\comma
$$
and each of these summands equals $-1$.
If $I_u(t)<I_u(0)$, the non-zero summands are
those with labels
$$
 i=I_u(t)+1,\ldots,I_u(0)\comma
$$
and each of these summands equals $1$.
If the two indices coincide, every summand is zero.
Thus, in all cases, only finitely many summands are non-zero,
and their sum is $I_u(0)-I_u(t)=J(u,t)$.
\end{proof}

\subsection{Compact linear statistics and the dynamic structure factor}
For a real-valued $f\in\mcS(\R)$, define the centred linear statistic
\begin{align}\label{e:xi}
\mcN_t(f)\defi\int_{\R}f(x)\bigl(\xi_t(\diff x)-\diff x\bigr)\comma
\qquad \xi_t=\sum_i\delta_{X_i(t)}\fstop
\end{align}
As the intensity of $\xi_t \sim \sine_2$ is $1$ for every~$t \ge 0$, we have  $\mbbE[\mcN_t(f)]=0$. For real random variables~$U$ and~$V$, we use the notation
$$
\Cov(U,V)\defi\mbbE[UV]-\mbbE[U]\mbbE[V]\comma
\qquad \Var(U)\defi\Cov(U,U)\fstop
$$

\begin{prop}
\label{prop:cov-FT}\label{prop:S-formula}
For each $t\ge0$, there is an even, real-valued function $S(k,t)$,
called the \emph{dynamic structure factor}, such that, for all
real-valued $f,g\in\mcS(\R)$,
\begin{equation}\label{eq:cov-FT}
\Cov(\mcN_t(f),\mcN_0(g))
=\frac{1}{2\pi}\int_{\R}
 \widehat f(k)\widehat g(-k)S(k,t)\diff k\fstop
\end{equation}
In particular,
\begin{equation}\label{eq:inc-var}
\Var(\mcN_t(f)-\mcN_0(f))
=\frac{1}{\pi}\int_{\R}|\widehat f(k)|^2
 \bigl(S(k,0)-S(k,t)\bigr)\diff k\fstop
\end{equation}
For $k\ge0$ and $t\ge0$, this function is given explicitly by
\begin{equation}\label{eq:S-Fermi}
S(k,t)=\frac{1}{2\pi}\int_{-\pi}^{\pi}
 \1_{\{|q+k|>\pi\}}
 \exp\!\left(-\frac{t}{2}\bigl((q+k)^2-q^2\bigr)\right)\diff q\fstop
\end{equation}
For $t>0$ and $k>0$, evaluating the integral gives
\begin{equation}\label{eq:S-closed}
S(k,t)=
\begin{cases}
\dfrac{1}{2\pi\,t\,k}
 \Bigl(e^{-t k(\pi-\frac{k}{2})}-e^{-t k(\pi+\frac{k}{2})}\Bigr),
 &0<k<2\pi,\\[3mm]
\dfrac{1}{2\pi\,t\,k}
 \Bigl(e^{-t k(\frac{k}{2}-\pi)}-e^{-t k(\frac{k}{2}+\pi)}\Bigr),
 &k\ge2\pi.
\end{cases}
\end{equation}
At $t=0$,
\begin{equation}\label{eq:S0}
S(k,0)=
\begin{cases}
\dfrac{k}{2\pi},&0\le k\le2\pi,\\
1,&k\ge2\pi,
\end{cases}
\qquad S(-k,0)=S(k,0)\fstop
\end{equation}
\end{prop}

\begin{proof}
We first express the covariance in terms of  spatial variables, then compute
its Fourier transform. A Fourier-inversion argument then
concludes~\eqref{eq:cov-FT} at all times, including $t=0$.
Finally, we evaluate the integral defining $S$.

\medskip
\noindent\emph{Step 1: the spatial covariance and the equal-time correction.}
Fix $t>0$. By expanding the two-by-two determinant
in~\eqref{eq:extended-sine-correlations} and using the unit
intensity property at both times, we obtain
$$
\rho^{(2)}_{\tr}((t,x),(0,y))
\defi\rho^{(2)}((t,x),(0,y))-1
=-K_{\sin}(t,x;0,y)K_{\sin}(0,y;t,x)\fstop
$$
Since the two observations lie on different time layers, the
definition of the two-point correlation density provides
$$
\mbbE\left[
 \left(\int_\R f(x)\xi_t(\diff x)\right)
 \left(\int_\R g(y)\xi_0(\diff y)\right)
\right]
=\iint_{\R^2}f(x)g(y)
 \rho^{(2)}((t,x),(0,y))\diff x\diff y\fstop
$$
Note that $\mbbE[\int_\R f(x)\xi_t(\diff x)]=\int_\R f(x)\diff x$ and $\mbbE[\int_\R g(x)\xi_0(\diff x)]=\int_\R g(y)\diff y$ by the stationarity~$\xi_t\sim\sine_2$ and the unit intensity.
By subtracting their product and using
the spatial translation invariance of the kernel, we have
$$
\begin{aligned} \label{eq:CFF}
\Cov(\mcN_t(f),\mcN_0(g))
&=\iint_{\R^2}f(x)g(y)
 \bigl(\rho^{(2)}((t,x),(0,y))-1\bigr)\diff x\diff y\\
&=\iint_{\R^2}f(x)g(y)\mcC_t(x-y)\diff x\diff y\comma
\end{aligned}
$$
where
\begin{equation}\label{e:Ct}
\mcC_t(z)\defi-K_{\sin}(t,z;0,0)K_{\sin}(0,0;t,z)\fstop
\end{equation}

For fixed $t>0$, the cosine formula~\eqref{eq:extended-sine-kernel}
shows that both kernel factors are real, even and bounded in $z$.
We now control the tails. For so doing,  we take integration by parts, which gives for $z \ge 0$
$$
\begin{aligned}
K_{\sin}(0,0;\,t,z)
&=\frac{e^{t\pi^2/2}\sin(\pi z)}{\pi z}
  -\frac{t}{\pi z}\int_0^\pi
       k e^{tk^2/2}\sin(kz)\diff k\comma\\
K_{\sin}(t,z;\,0,0)
&=\frac{e^{-t\pi^2/2}\sin(\pi z)}{\pi z}
  -\frac{t}{\pi z}\int_\pi^\infty
       k e^{-tk^2/2}\sin(kz)\diff k\fstop
\end{aligned}
$$
The boundary term at $k=0$ vanishes because $\sin(0)=0$ and 
the boundary term at infinity vanishes because $t>0$ and the
Gaussian weight tends to zero. Moreover, we have the following identities:
$$
 t\int_0^\pi k e^{tk^2/2}\diff k=e^{t\pi^2/2}-1\comma
 \qquad
 t\int_\pi^\infty k e^{-tk^2/2}\diff k=e^{-t\pi^2/2}\fstop
$$
Using the bound~$|\sin(kz)|\le1$, we obtain a bound $C(t)/|z|$ for each
kernel. Their negative product $\mcC_t$ is  real
and even, bounded near zero, and $O_t(|z|^{-2})$ at infinity.
In particular, $\mcC_t\in L^1(\R)$.

Furthermore, the covariance at $t=0$ is given as 
$$
\Cov(\mcN_0(f),\mcN_0(g))
=\int_\R f(x)g(x)\diff x
 -\iint_{\R^2}f(x)g(y)K_{\sin}(x-y)^2\diff x\diff y\fstop
$$
Thus the covariance measure in the variable $z=x-y$ is
$$
\mcC_0(\diff z)\defi\delta_0(\diff z)-K_{\sin}(z)^2\diff z\fstop
$$
We write
$\mcC_t(\diff z)=\mcC_t(z)\diff z$ when $t>0$. Note that they are finite signed measures. Indeed,  for $t>0$ this follows from
the integrability just proved above, and at $t=0$ it follows from the integrability~$K_{\sin}\in L^2(\R)$.

\medskip
\noindent\emph{Step 2: computation of the Fourier transform.}
For every $t\ge0$, define
$$
 S(k,t)\defi\int_\R e^{-ikz}\mcC_t(\diff z)\fstop
$$
By definition, $S(k,t)=\widehat\mcC_t(k)$ for $t>0$.
Since the covariance measures are real and even, so is $S$.
 Using~\eqref{eq:ext-sine-FT}, \eqref{eq:kappa-KT} and~\eqref{e:Ct}  we have 
$$
\begin{aligned}
S(k,t) = \widehat\mcC_t(k)
&=-\frac{1}{2\pi}\int_\R
  \kappa_{t,0}(k-q)\kappa_{0,t}(-q)\diff q\\
&=-\frac{1}{2\pi}\int_\R
  \kappa_{t,0}(k+q)\kappa_{0,t}(q)\diff q\\
&=\frac{1}{2\pi}\int_\R
  \1_{\{|q|\le\pi\}}\1_{\{|q+k|>\pi\}}
  \exp\!\left(-\frac{t}{2}\bigl((q+k)^2-q^2\bigr)\right)\diff q\fstop
\end{aligned}
$$
The first integral follows from the product-convolution formula.
The remaining equalities follow from the change of variable
$q\mapsto-q$ and substitution of the Fourier multipliers.
This  proves~\eqref{eq:S-Fermi}
for $t>0$.

At $t=0$, the diagonal measure gives
$S(k,0)=1-\widehat{K_{\sin}^2}(k)$.
Since $\widehat K_{\sin}(q)=\1_{\{|q|\le\pi\}}$, the same convolution formula yields
$$
\begin{aligned}
S(k,0)
&=1-\frac{1}{2\pi}\int_{-\pi}^{\pi}
       \1_{\{|k-q|\le\pi\}}\diff q
&=1-\frac{1}{2\pi}\int_{-\pi}^{\pi}
       \1_{\{|k+q|\le\pi\}}\diff q
&=\frac{1}{2\pi}\int_{-\pi}^{\pi}
       \1_{\{|q+k|>\pi\}}\diff q\fstop
\end{aligned}
$$
Here we again changed $q$ to $-q$, then used
$1=(2\pi)^{-1}\int_{-\pi}^{\pi}\diff q$.
This proves~\eqref{eq:S-Fermi} also at $t=0$.

\medskip
\noindent\emph{Step 3: the covariance formula.}
Let $g^\vee(x)\defi g(-x)$ and set $h=f*g^\vee$.
For every $t\ge0$,
$$
 \Cov(\mcN_t(f),\mcN_0(g))
 =\int_\R h(z)\mcC_t(\diff z)\fstop
$$
At $t=0$, the atom at zero gives
$h(0)=\int_\R f(x)g(x)\diff x$ as required.
Since $h$ is a Schwartz function and $\mcC_t$ is a finite
signed measure, by the Fourier inversion formula and Fubini's theorem, we obtain
$$
\begin{aligned}
\Cov(\mcN_t(f),\mcN_0(g))
&=\frac{1}{2\pi}\int_\R
       \widehat h(k)S(-k,t)\diff k
&=\frac{1}{2\pi}\int_\R
       \widehat f(k)\widehat g(-k)S(k,t)\diff k\fstop
\end{aligned}
$$
We used $\widehat{g^\vee}(k)=\widehat g(-k)$ and the evenness
of $S$. This proves~\eqref{eq:cov-FT} simultaneously for
distinct and equal time-points.

By the stationarity,
$$
\Var(\mcN_t(f)-\mcN_0(f))
=2\Var(\mcN_0(f))-2\Cov(\mcN_t(f),\mcN_0(f))\fstop
$$
Applying~\eqref{eq:cov-FT} and using
$\widehat f(k)\widehat f(-k)=|\widehat f(k)|^2$ for real-valued functions~$f$, we proved~\eqref{eq:inc-var}.

\medskip
\noindent\emph{Step 4: evaluation of the structure factor.}
It remains to evaluate~\eqref{eq:S-Fermi} for $k\ge0$.
If $k=0$,  $S(0,t)=0$.
For $k>0$, write
$$
 \frac{(q+k)^2-q^2}{2}=kq+\frac{k^2}{2}\fstop
$$
Suppose first that $0<k<2\pi$. For $q\in[-\pi,\pi]$, we have
$q+k>-\pi$, so $|q+k|>\pi$ is equivalent to $q>\pi-k$.
Consequently, for $t>0$,
$$
\begin{aligned}
S(k,t)
&=\frac{e^{-tk^2/2}}{2\pi}
       \int_{\pi-k}^{\pi}e^{-tkq}\diff q
&=\frac{e^{-tk^2/2}}{2\pi tk}
       \Bigl(e^{-tk(\pi-k)}-e^{-tk\pi}\Bigr)
&=\frac{1}{2\pi tk}
       \Bigl(e^{-tk(\pi-k/2)}-e^{-tk(\pi+k/2)}\Bigr)\fstop
\end{aligned}
$$
This is the first line of~\eqref{eq:S-closed}.

If $k\ge2\pi$, then $q+k\ge\pi$ throughout $[-\pi,\pi]$.
The equality occurs only at $q=-\pi$ when $k=2\pi$, which does
not affect the integral. Thus we have
% the integration domain is the
%whole interval, up to a set of measure zero, and
$$
\begin{aligned}
S(k,t)
&=\frac{e^{-tk^2/2}}{2\pi}
       \int_{-\pi}^{\pi}e^{-tkq}\diff q
&=\frac{e^{-tk^2/2}}{2\pi tk}
       \Bigl(e^{tk\pi}-e^{-tk\pi}\Bigr)
&=\frac{1}{2\pi tk}
       \Bigl(e^{-tk(k/2-\pi)}-e^{-tk(k/2+\pi)}\Bigr)\fstop
\end{aligned}
$$
This is the second line of~\eqref{eq:S-closed}.

Finally, at $t=0$ the exponential is  equal to one, so $S(k,0)$ is
the length of the integration domain divided by $2\pi$.
As the lengths  are $k$ for $0\le k\le2\pi$ and
$2\pi$ for $k\ge2\pi$, we conclude~\eqref{eq:S0}.
\end{proof}

\begin{lem}\label{lem:S-bounds}
For every $t\ge0$ and $k\in\R$,
\begin{enumerate}
\item $0\le S(k,t)\le S(k,0)\le1$;
\item $0\le S(k,0)-S(k,t)\le\min\{1,tk^2/2\}$.
\end{enumerate}
In particular, for $|k|\le1$ we have $S(k,0)-S(k,t)\le C(t)k^2$ with $C(t)=t/2$.
\end{lem}

\begin{proof}
By the evenness, it suffices to consider $k\ge0$. On the integration set in~\eqref{eq:S-Fermi}, we have $|q|\le\pi<|q+k|$, so
$$
E_k(q)\defi\frac{(q+k)^2-q^2}{2}=qk+\frac{k^2}{2}\ge0\fstop
$$
Consequently, the exponential in~\eqref{eq:S-Fermi} lies in $[0,1]$. This proves the first assertion and the bound $S(k,0)-S(k,t)\le1$.

For the remaining bound, use $1-e^{-v}\le v$ for $v\ge0$ to obtain
$$
S(k,0)-S(k,t)
\le\frac{t}{2\pi}\int_{-\pi}^{\pi}\1_{\{|q+k|>\pi\}}E_k(q)\diff q\fstop
$$
If $0\le k\le2\pi$, then
$$
\int_{\pi-k}^{\pi}\left(qk+\frac{k^2}{2}\right)\diff q=\pi k^2\fstop
$$
If $k\ge2\pi$, then
$$
\int_{-\pi}^{\pi}\left(qk+\frac{k^2}{2}\right)\diff q=\pi k^2\fstop
$$
Thus the same upper bound $tk^2/2$ holds in both cases.
\end{proof}

For later use, we note here  the time-dependent factor in the low-frequency formula. For $t>0$ and $0<|k|<2\pi$, define the even function
\begin{align}\label{d:B}
B(k,t)\defi e^{-t|k|(\pi-\frac{|k|}{2})}\,
\frac{1-e^{-tk^2}}{tk^2}\comma \qquad
S(k,t)=\frac{|k|}{2\pi}B(k,t)\fstop
\end{align}
The two factors in this definition belong to $[0,1]$. Moreover, Lemma~\ref{lem:S-bounds} gives
$$
0\le1-B(k,t)
=\frac{2\pi}{|k|}\bigl(S(k,0)-S(k,t)\bigr)
\le\pi t|k|\fstop
$$

\begin{lem}\label{lem:B-bounds}
Let $B(k,t)$ be as in~\eqref{d:B}, and fix $k_0\in(0,2\pi)$. There exists $c=c(k_0)>0$ such that, for every $t>0$ and $0<|k|\le k_0$,
$$
0\le B(k,t)\le e^{-ct|k|}\fstop
$$
\end{lem}

\begin{proof}
The inequality $1-e^{-v}\le v$ implies
$$
0\le\frac{1-e^{-tk^2}}{tk^2}\le1\fstop
$$
Set $c\defi\pi-k_0/2>0$. Since $\pi-|k|/2\ge c$, by applying these bounds to~\eqref{d:B}, we obtain 
$$
0\le B(k,t)\le e^{-t|k|(\pi-|k|/2)}\le e^{-ct|k|}\comma
$$
as required.
\end{proof}

\begin{lem}[Extension to compactly supported measurable tests]\label{lem:compact-covariance}
The statistic in~\eqref{e:xi} is also well-defined and belongs to $L^2$
for every bounded, real-valued, compactly supported measurable function
$f$. Equations~\eqref{eq:cov-FT} and~\eqref{eq:inc-var} remain valid
for such test functions.
\end{lem}

\begin{proof}
For a bounded Borel set $B$ and a time~$s$, the one-/two-point correlations of the sine point process give
\begin{equation}\notag %\label{eq:local-count-second-moment}
\mbbE[\xi_s(B)^2]
=|B|+\int_{B^2}\bigl(1-|K_{\sin}(x-y)|^2\bigr)\diff x\diff y
\le |B|+|B|^2\fstop
\end{equation}
%Thus a linear statistic supported in $B$ has a finite second moment.
For a real-valued function~$h$ supported in $B$, 
$$
\begin{aligned}
\mbbE\left|\int_{\R}h(x)\xi_s(\diff x)\right|^2
&\le\int_B|h(x)|^2\diff x
  +\int_{B^2}|h(x)h(y)|\diff x\diff y\\
&=\|h\|_{L^2}^2+\|h\|_{L^1}^2
\le(1+|B|)\|h\|_{L^2}^2\fstop
\end{aligned}
$$
Here we used $0\le\rho^{(2)}\le1$ and the Cauchy--Schwarz inequality on~$B$.

We now choose smooth compactly supported functions $f_n$ converging to $f$
in $L^2$ with all of their supports contained in one bounded interval $B$.
 The preceding estimate shows that
$\mcN_s(f_n)\to\mcN_s(f)$ in $L^2$ at each observation time.
For two test functions, the Cauchy--Schwarz inequality therefore implies the convergence
of the covariances on the left-hand side of~\eqref{eq:cov-FT}.
On the righthand side, the Plancherel formula and the bound~$|S(k,t)|\le1$ show that the
bilinear integral is continuous in the $L^2$ norms.
Passing to the limit proves~\eqref{eq:cov-FT}. Its increment form~\(\mathcal N_t(f)-\mathcal N_0(f)\)
 gives~\eqref{eq:inc-var}.
\end{proof}

\subsection{Weighted currents and the distant boundary}
Fix $\chi\in C_c^\infty(\R)$ with $0\le\chi\le1$,
$\chi=1$ on $[-1,1]$ and $\chi=0$ outside $[-2,2]$.
For $L\ge1$, we define $\chi_L(x)\defi\chi(x/L)$ and
\begin{equation}\label{eq:boundary-weight}
\psi(y)\defi\chi'(y)\1_{\{y>0\}}\fstop
\end{equation}
The multiplication of the indicator function does not causes loss of smoothness: $\chi'$ is zero in a
neighbourhood of zero, so $\psi\in C_c^\infty(\R)$ and its support
is contained in $[1,2]$.

We first prove a general estimate for  weighted currents.
This estimate will control the boundary correction of the half-line
cut-off.

\begin{lem}[Variance formula]\label{l:WCV}
Fix $t>0$. Let $a\in C_c^\infty(\R)$ be real-valued and set
$$
A(x)\defi \int_{-\infty}^x a(u)\diff u\comma \qquad
W_a(t)\defi \int_{\R}a(u)J(u,t)\diff u\fstop
$$
Then $W_a(t)\in L^2$, $\mbbE[W_a(t)]=0$, and
\begin{equation}\label{eq:weighted-current-var}
\Var(W_a(t))
=\frac1\pi\int_{\R}\frac{|\widehat a(k)|^2}{k^2}
\big(S(k,0)-S(k,t)\big)\diff k\fstop
\end{equation}
The integrand is bounded near $k=0$, and its value at $k=0$ does not affect the integral.
Moreover, with the cut-off $\chi_L$ fixed above,
$$
\mcN_t(\chi_LA)-\mcN_0(\chi_LA)
\xrightarrow[L\to\infty]{L^2}W_a(t)\fstop
$$
\end{lem}

\begin{proof}
\noindent\textit{Step 1: construction of the $L^2$ limit.}
Set~$c\defi\int_{\R}a(u)\diff u$, $A_L\defi\chi_LA$, and define
$$
X_{a,L}\defi\mcN_t(A_L)-\mcN_0(A_L)\fstop
$$
We take $L$ so large that $\supp\ a\subset[-L,L]$.
On the negative support of $\chi_L'$, the antiderivative $A$ is zero, and 
on its positive support, $A$ is equal to $c$. By differentiating the product
$A_L=\chi_LA$, we therefore have
$$
A_L'=a+cL^{-1}\psi(\cdot/L)\fstop
$$
By integration by parts and a change of variables,  we have that for $k\ne0$,
$$
\widehat A_L(k)
=\frac{\widehat a(k)+c\widehat\psi(Lk)}{ik}\fstop
$$
The Riemann--Lebesgue lemma implies $\widehat\psi(Lk)\to0$ for
fixed $k\ne0$. Also,
$$
|\widehat A_L(k)|
\le\frac{\|a\|_{L^1}+|c|\|\psi\|_{L^1}}{|k|}\fstop
$$
Applying~\eqref{eq:inc-var} to $A_L-A_{L'}$ shows that
$$
\mbbE|X_{a,L}-X_{a,L'}|^2
=\frac1\pi\int_{\R}
|\widehat A_L(k)-\widehat A_{L'}(k)|^2
\bigl(S(k,0)-S(k,t)\bigr)\diff k\fstop
$$
The integrand tends to zero for every $k\ne0$ as $L,L'\to\infty$.
By Lemma~\ref{lem:S-bounds}, it is bounded by a constant times
$\min\{t/2,k^{-2}\}$, which is integrable for fixed $t$.
By the dominated convergence,  $(X_{a,L})_L$ is Cauchy in $L^2$.
We denote its limit by $\widetilde W_a(t)$.
As $X_{a,L}$ are centred, the limit is centred as well.
The same dominating function allows us to passage to the limit in their
variances and, therefore, \eqref{eq:weighted-current-var} follows with
$\widetilde W_a(t)$ in place of $W_a(t)$.
The linearity of the compact statistics also gives linearity of the function~$a\mapsto\widetilde W_a(t)$.

\smallskip
\noindent\textit{Step 2: reduction to a single discrepancy.}
Let $ I_\alpha $ be as in \eqref{:22}. 
Choose $\alpha<\beta$ with $\supp~a\subset[\alpha,\beta]$,
and set
$$
m\defi\min\{I_\alpha(0),I_\alpha(t)\}\comma\qquad
n\defi\max\{I_\beta(0),I_\beta(t)\}\fstop
$$
For $i\le m$, both endpoint positions are at or below $\alpha$;
for $i>n$, both are above $\beta$. The differences of indicator functions in
Lemma~\ref{lem:finite-sum} therefore vanish outside the fixed finite
set $\{m+1,\ldots,n\}$, uniformly for $u\in[\alpha,\beta]$.
We may integrate that finite sum against $a$ and obtain
\begin{equation}\label{eq:weighted-current-pathwise}
W_a(t)=\sum_i\bigl(A(X_i(t))-A(X_i(0))\bigr)\fstop
\end{equation}
In particular, $W_a(t)$ is finite almost surely. 

Now we  identify $W_a(t)$, defined pathwise, 
with the $L^2$-limit~$\widetilde W_a(t)$ constructed in Step~1.
If $\int _{\R } a du =0$, then $A$ is compactly supported. For all sufficiently
large $L$, we therefore have $\chi_LA=A$, and the finite sum above
is  $X_{a,L}$. Passing to the $L^2$ limit gives
$W_a(t)=\widetilde W_a(t)$ almost surely.

Choose a real-valued function~$a_0\in C_c^\infty(\R)$ with $\int _{\R} a_0 du =1$ and set
$D(t)\defi W_{a_0}(t)-\widetilde W_{a_0}(t)$.
Recall that $c=\int_\R a \diff u$. Since $\int_\R a_0 \diff u=1$,
the function $a-ca_0$ has zero integral.
The zero-integral case proved above therefore gives
$W_{a-ca_0}(t)=\widetilde W_{a-ca_0}(t)$ almost surely.
By the linearity of both maps, this equality becomes
$$
 W_a(t)-cW_{a_0}(t)
 =\widetilde W_a(t)-c\widetilde W_{a_0}(t)\fstop
$$
By the definition of $D(t)$, we have
$$
 W_a(t)-\widetilde W_a(t)
 =c\bigl(W_{a_0}(t)-\widetilde W_{a_0}(t)\bigr)
 =\left(\int_\R a(u)\diff u\right)D(t)\fstop
$$
We shall show that $D(t)=0$.

\smallskip
\noindent\textit{Step 3: spatial ergodicity of the path law.}
To show in Step~4 that $D(t)$ is deterministic, we here prove the spatial
ergodicity of the unlabelled path law on $[0,t]$.
Let $\theta_h$ translate every particle position by $h$.  The path law
is translation invariant by Section~\ref{subsec:extended-sine-current}.
We prove the stronger mixing identity
\begin{align} \label{e:MID}
 \mbbE[F(G\circ\theta_n)]
 \longrightarrow \mbbE[F]\mbbE[G]\comma
 \qquad n\to\infty\comma
\end{align}
for all bounded measurable functions $F,G$ of the unlabelled paths, i.e., the space of continuous paths in the configuration space over $\R$.

We adapt the argument in~\cite[proof of Theorem~7]{Soshnikov2000DRPF}
to the finite-time field in~\eqref{eq:extended-sine-correlations} with our non-Hermitian kernel.
We first take $F=\Phi(U)$ and $G=\Psi(V)$, where $\Phi,\Psi$ are
bounded functions and
$$
 U\defi\bigl(\xi_{s_i}(I_i)\bigr)_{i=1}^p\comma\qquad
 V\defi\bigl(\xi_{r_j}(J_j)\bigr)_{j=1}^q\comma\qquad
 V_n\defi V\circ\theta_n
 =\bigl(\xi_{r_j}(J_j-n)\bigr)_{j=1}^q\fstop
$$
Here $p,q\ge1$, the times $s_i,r_j\in[0,t]$ are fixed,
and $I_i,J_j$ are bounded Borel sets; $\xi_s(I)$ counts the
particles in $I$ at time $s$. By the translation invariance, we have
$\law(V_n)=\law(V)$.

\smallskip
\noindent\emph{(i) Factorisation of count moments.}
For fixed observation times $s$ and $r$, the multiplier
in~\eqref{eq:kappa-KT} is integrable. By the Riemann--Lebesgue
lemma, we have that, for bounded intervals $I,J$,
$$
 \sup_{\substack{x\in I\\y\in J}}
 \Bigl(
 |K_{\sin}(s,x;\,r,y-n)|
 +|K_{\sin}(r,y-n;\,s,x)|
 \Bigr)\longrightarrow0\fstop
$$
Choose $I,J$ containing all $I_i,J_j$, respectively.
For integers $k,\ell\ge1$, consider $z_i=(\sigma_i,x_i)$, $1\le i\le k$, and
$w_j^{(n)}=(\tau_j,y_j-n)$, $1\le j\le\ell$, with fixed times
$\sigma_i,\tau_j\in[0,t]$ and $x_i\in I$, $y_j\in J$.
The  correlation matrix in the order
$z_1,\ldots,z_k$ and $w_1^{(n)},\ldots,w_\ell^{(n)}$ has the following blocks
$$
 \begin{aligned}
 \mathsf A_{ii'}
 &=K_{\sin}(\sigma_i,x_i;\,\sigma_{i'},x_{i'})\comma
 &\mathsf D_{jj'}
 &=K_{\sin}(\tau_j,y_j;\,\tau_{j'},y_{j'}) \comma\\
 (\mathsf U_n)_{ij}
 &=K_{\sin}(\sigma_i,x_i;\,\tau_j,y_j-n)\comma
 &(\mathsf V_n)_{ji}
 &=K_{\sin}(\tau_j,y_j-n;\,\sigma_i,x_i)\fstop
 \end{aligned}
$$
By~\eqref{eq:extended-sine-correlations}, the joint correlation
density is explicitly given as
$$
 \rho^{(k+\ell)}(z_1,\ldots,z_k,w_1^{(n)},\ldots,w_\ell^{(n)})
 =\det\begin{pmatrix}
  \mathsf A&\mathsf U_n\\
  \mathsf V_n&\mathsf D
 \end{pmatrix}\fstop
$$
In $\mathsf D$, the translation invariance cancels the common
shift $-n$. Thus $\mathsf A,\mathsf D$ are independent of $n$,
while every entry of $\mathsf U_n,\mathsf V_n$ tends uniformly
to zero by the preceding decay estimate. Hence
$$
 \det\begin{pmatrix}
  \mathsf A&\mathsf U_n\\
  \mathsf V_n&\mathsf D
 \end{pmatrix}
 \longrightarrow
 \det\begin{pmatrix}
  \mathsf A&0\\
  0&\mathsf D
 \end{pmatrix}
 =\det(\mathsf A)\det(\mathsf D)\fstop
$$
The convergence is uniform over the spatial coordinates.
Integrating the both sides, we obtain factorisation of mixed factorial moments of $(U,V_n)$. 

\smallskip
\noindent\emph{(ii) From count moments to local mixing.}
By expanding the moments, we obtain
$$
 \mbbE e^{\lambda\xi_s(B)}
 \le\exp\!\bigl(|B|(e^\lambda-1)\bigr)\comma
 \qquad\lambda\ge0\fstop
$$
For any $m$ local counts $N_j=\xi_{u_j}(B_j)$, by H\"older's
inequality with all exponents equal to $m$, we have 
$$
 \mbbE e^{\lambda\sum_{j=1}^mN_j}
 \le\prod_{j=1}^m
       \bigl(\mbbE e^{m\lambda N_j}\bigr)^{1/m}
 \le\exp\!\left(
       \frac{e^{m\lambda}-1}{m}\sum_{j=1}^m|B_j|
       \right)\fstop
$$
This bound is unchanged by translating any $B_j$, so it applies
uniformly to $(U,V_n)$. It justifies  limits in their
joint moment generating functions near zero.
The moment factorisation in~(i) then identifies the limiting law
as $\law(U)\otimes\law(V)$.
Since these vectors are integer-valued, convergence of their laws
also gives convergence of bounded expectations. This proves the
displayed mixing identity~\eqref{e:MID} for $F=\Phi(U)$ and $G=\Psi(V)$.

\smallskip
\noindent\emph{(iii) Extension to the path law and ergodicity.}
The \emph{path sigma-field} is the \(\sigma\)-field on~$\Omega$ generated by the unlabelled process up to time \(t\), viz.,
$$
 \mathcal G_t^\xi\defi
 \sigma\bigl(\xi_s(B):0\le s\le t,\;
 B\subset\R\ \text{bounded Borel}\bigr)\fstop
$$
It consists of the events determined by the unlabelled configurations
observed throughout $[0,t]$.
Put $\mathcal D_t=(\mathbb Q\cap[0,t])\cup\{t\}$.
By the continuity of the unlabelled paths, we have that $\mathcal G_t^\xi$
is generated by the countable family
$\xi_s((a,b])$ with $s\in\mathcal D_t$ and
$a,b\in\mathbb Q$, $a<b$.
Thus bounded $\mathcal G_t^\xi$-measurable functions can be approximated
in $L^2$ by bounded functions of finitely many such counts.
Cauchy--Schwarz and translation invariance make the approximation
errors in the mixing identity uniform in $n$.
It therefore holds for all bounded measurable path functions.
In particular, if an event $E$ is invariant under $\theta_1$ (in particular, it is invariant under~$\theta_n$ for every $n \in \mathbb N$), then
$$
 \mbbP(E)
 =\mbbE\bigl[\1_E(\1_E\circ\theta_n)\bigr]
 \longrightarrow\mbbP(E)^2\fstop
$$
Hence $\mbbP(E)\in\{0,1\}$, proving the spatial ergodicity
needed in Step~4.

\smallskip
\noindent\textit{Step 4: identification of the current.}
The ordered non-colliding trajectories can be recovered from the
unlabelled path up to an overall integer relabelling, which does not change 
the current. Thus $W_a(t)$ is $\mathcal G_t^\xi$-measurable, as is
$\widetilde W_a(t)$ by its construction as an $L^2$ limit.
Consequently, $D(t)$ is measurable with respect to the completion
of $\mathcal G_t^\xi$. 
For $h\in\R$, set~$a_h(u)\defi a(u+h)$.
The translation of paths gives
$$
W_a(t)\circ\theta_h=W_{a_h}(t)\fstop
$$
The same identity holds for the $L^2$-limit $\widetilde W_a(t)$.
Indeed, recalling and setting \[
\begin{aligned}
X_{a,L}\circ\theta_h-X_{a_h,L}
&=\mathcal N_t(d_{L,h})-\mathcal N_0(d_{L,h}) \comma \quad 
d_{L,h}(x):=A(x+h)\bigl(\chi_L(x+h)-\chi_L(x)\bigr).
\end{aligned}
\] 
The fundamental theorem of calculus gives
$\chi_L(x+h)-\chi_L(x)=\int_0^h\chi_L'(x+r)\diff r$.
Since $\|\chi_L'\|_{L^2}=L^{-1/2}\|\chi'\|_{L^2}$,
Minkowski's inequality and the translation invariance of the $L^2$ norm imply
$$
\|d_{L,h}\|_{L^2}
\le\|A\|_\infty |h|L^{-1/2}\|\chi'\|_{L^2}\fstop
$$
Using~\eqref{eq:inc-var}, $S(k,0)-S(k,t)\le1$ and the Plancherel formula,
we obtain
$$
\mbbE\big|
\mcN_t(d_{L,h})-\mcN_0(d_{L,h})
\big|^2
\le2\|d_{L,h}\|_{L^2}^2\xrightarrow{L \to \infty}0\fstop
$$
The translation preserves $L^2$ convergence, so
$\widetilde W_a(t)\circ\theta_h=\widetilde W_{a_h}(t)$.

Since $\int(a_0)_h=1$, applying Step~2, we have 
$$
D(t)\circ\theta_h
=W_{(a_0)_h}(t)-\widetilde W_{(a_0)_h}(t)=D(t)\fstop
$$
Thus $D(t)$ is invariant under $\theta_1$. For every rational $r$,
the event $\{D(t)\le r\}$ is invariant under $\theta_1$, thus the spatial ergodicity implies that such an event has  probability zero or
one for every $r$. Since $D(t)$ is finite almost surely, these events show that
it is a deterministic constant,

For fixed $t>0$, define the time-reversal map on the path space by
$$
(\mcR_t\omega)(s)\defi\omega(t-s)\comma\qquad 0\le s\le t\fstop
$$
It follows the same trajectory in the backward way and interchange the initial
and terminal states. Following the particle labels backwards therefore
reverses the current while exchanging the endpoint statistics reverses
their difference. Namely, 
$$
W_a(t)\circ\mcR_t=-W_a(t)\comma\qquad
X_{a,L}\circ\mcR_t=-X_{a,L}\fstop
$$

By the reversibility,  $\mcR_t$ preserves the path law, so the composition
with $\mcR_t$ preserves the $L^2$ norm (see, e.g.,~\cite[Lemma 5.7.1]{FOT11}). Hence
$$
\|X_{a,L}\circ\mcR_t-\widetilde W_a(t)\circ\mcR_t\|_{L^2}
=\|X_{a,L}-\widetilde W_a(t)\|_{L^2}
\longrightarrow0\fstop
$$
On the other hand, $X_{a,L}\circ\mcR_t=-X_{a,L}$ converges in $L^2$
to $-\widetilde W_a(t)$. The uniqueness of the $L^2$ limit therefore gives
$$
\widetilde W_a(t)\circ\mcR_t=-\widetilde W_a(t)
\qquad\text{almost surely}\fstop
$$
Recalling $D(t)=W_{a_0}(t)-\widetilde W_{a_0}(t)$ and the corresponding
sign identity for $W_{a_0}(t)$, we obtain
\begin{align}\label{eZR}
D(t)\circ\mcR_t
=-W_{a_0}(t)+\widetilde W_{a_0}(t)
=-D(t)
\qquad\text{almost surely}\fstop
\end{align}
Noting that $D(t)=D(t)\circ\mcR_t$ as  $D(t)$ is a constant almost surely, \eqref{eZR} implies  $D(t)=0$.
Hence $W_a(t)=\widetilde W_a(t)$, and all assertions follow from
Step~1.
\end{proof}

\subsection{Approximation and square-integrability of the current}
Let $H(x)\defi\1_{\{x>0\}}$ and, for $R\ge1$, recalling that $\chi_R$ is the cut-off function defined in~\eqref{eq:boundary-weight}, set
\begin{equation}\label{eq:sharp-current-cutoff}
f_R(x)\defi\chi_R(x)H(x)\comma\qquad
Z_R(t)\defi\mcN_t(f_R)-\mcN_0(f_R)\fstop
\end{equation} 
Because $f_R$ is bounded and compactly supported,
Lemma~\ref{lem:compact-covariance} shows that $Z_R(t)\in L^2$ and
$\mbbE[Z_R(t)]=0$.

\begin{prop}\label{prop:direct-current-L2}
For every fixed $t>0$, $J(t)\in L^2$ and $\mbbE[J(t)]=0$.
Moreover, for every $R\ge1$,
\begin{equation}\label{eq:sharp-current-error}
\mbbE|Z_R(t)-J(t)|^2
\le\frac{t}{R}\|\psi\|_{L^2}^2\comma
\end{equation}
where $\psi$ was defined in~\eqref{eq:boundary-weight}. 
In particular $Z_R(t)\to J(t)$ in $L^2$ as $R\to\infty$.
\end{prop}

\begin{proof}
\noindent\emph{Step 1: identify the boundary correction pathwise.}
Let $ \psi(y) = \chi'(y)\1_{\{y>0\}}$ be as in \eqref{eq:boundary-weight}. 
Define
$$
a_R(u)\defi\chi_R'(u)H(u)=R^{-1}\psi(u/R)\fstop
$$
This is a smooth compactly supported function: the derivative of the cut-off vanishes near the origin.
Its antiderivative is
$$
D_R(x)\defi\int_{-\infty}^x a_R(u)\diff u
=f_R(x)-H(x)\fstop
$$
Indeed, both sides are zero for $x\le0$, and for $x>0$, the integral is
$\chi_R(x)-\chi_R(0)=\chi_R(x)-1$.

Since $f_R$ is compactly supported, its endpoint statistic can be
written as a finite sum over labels of particles. The  summation for $J(t)$ is given by Lemma~\ref{lem:finite-sum} while the summation for $W_{a_R}(t)$ is given by
\eqref{eq:weighted-current-pathwise}.
By subtracting the first two sums and using the relation $D_R=f_R-H$, we obtain 
\begin{equation}\label{eq:sharp-current-boundary}
Z_R(t)-J(t)=W_{a_R}(t)\fstop
\end{equation}

\smallskip
\noindent\emph{Step 2: establish square integrability.}
By Lemma~\ref{l:WCV}, we have that $W_{a_R}(t)\in L^2$ and
$\mbbE[W_{a_R}(t)]=0$. For any fixed $R$,
the formula~\eqref{eq:sharp-current-boundary} expresses $J(t)$ as the difference
of the two centred random variables $Z_R(t)$ and $W_{a_R}(t)$, both of which belong to $L^2$
Thus $J(t)\in L^2$ and $\mbbE[J(t)]=0$.

\smallskip
\noindent\emph{Step 3: remove the distant boundary.}
We now apply~\eqref{eq:weighted-current-var} in Lemma~\ref{l:WCV} to $a_R$ and use the bound~$S(k,0)-S(k,t)\le tk^2/2$ from Lemma~\ref{lem:S-bounds}.
By the Plancherel formula and the change of variables $u=Rv$, we have
$$
\begin{aligned}
\mbbE|Z_R(t)-J(t)|^2
&=\frac1\pi\int_{\R}\frac{|\widehat a_R(k)|^2}{k^2}
       \bigl(S(k,0)-S(k,t)\bigr)\diff k
&\le\frac{t}{2\pi}\int_{\R}|\widehat a_R(k)|^2\diff k
=t\|a_R\|_{L^2}^2
=\frac{t}{R}\|\psi\|_{L^2}^2\fstop
\end{aligned}
$$
This proves~\eqref{eq:sharp-current-error} and the assertion about the  $L^2$
convergence.
\end{proof}

\begin{prop}\label{prop:Var-J}
For each fixed $t>0$,
\begin{equation}\label{eq:Var-J}
\Var(J(t))
=\frac{1}{\pi}\int_{\R}\frac{S(k,0)-S(k,t)}{k^2}\diff k\fstop
\end{equation}
\end{prop}

\begin{proof}
By Proposition~\ref{prop:direct-current-L2},
$\Var(Z_R(t))\to\Var(J(t))$.
For $k\ne0$, integration by parts on $[0,2R]$ gives
$$
\begin{aligned}
\widehat f_R(k)
&=\int_0^{2R}\chi_R(x)e^{-ikx}\diff x
&=\frac{1}{ik}+\frac{1}{ik}
  \int_0^{2R}\chi_R'(x)e^{-ikx}\diff x
=\frac{1+\widehat\psi(Rk)}{ik}\fstop
\end{aligned}
$$
The boundary contribution at zero is $1/(ik)$ because
$\chi_R(0)=1$; the contribution at $2R$ is zero.
The Riemann--Lebesgue lemma now gives
$\widehat f_R(k)\to1/(ik)$ for each fixed $k\ne0$, and
$$
|\widehat f_R(k)|\le\frac{1+\|\psi\|_{L^1}}{|k|}\fstop
$$
Consequently, by Lemma~\ref{lem:compact-covariance}, 
$$
\Var(Z_R(t))
=\frac1\pi\int_{\R}|\widehat f_R(k)|^2
       \bigl(S(k,0)-S(k,t)\bigr)\diff k\fstop
$$
By Lemma~\ref{lem:S-bounds}, the integrand is bounded, independently
of $R$, by
$$
(1+\|\psi\|_{L^1})^2
\min\left\{\frac t2,\frac1{k^2}\right\}\fstop
$$
For fixed $t$, this function is bounded near zero and integrable at
infinity. By the dominated convergence, we proved~\eqref{eq:Var-J}.
\end{proof}

\subsection{The logarithmic variance asymptotic}
%The exact variance formula is now established for the physical
%current. 
%The logarithm comes from frequencies between order $t^{-1}$
%and a fixed positive scale. Spohn~\cite{Spo87} established the leading logarithmic asymptotic. We give a direct proof for the pathwise current, including the bounded remainder in~\eqref{e:VF}.
\begin{thm}\label{thm:main}
As $t\to\infty$,
$$
\Var(J(t))=\frac{1}{\pi^2}\log t+O(1)\fstop
$$
\end{thm}

\begin{proof} In the following, we evaluate the right-hand side of the formula~\eqref{eq:Var-J}. 

\noindent\emph{Step 1: the small frequencies.}
Fix $k_0\in(0,1)$. By $0\le S(k,0)-S(k,t)\le1$, the part of
the integral in~\eqref{eq:Var-J} over $|k|>k_0$ is bounded by
$2/k_0$ independently of $t$. On $0<k\le k_0$, the factorisation
in~\eqref{d:B} implies
$$
\frac{S(k,0)-S(k,t)}{k^2}
=\frac{1}{2\pi}\frac{1-B(k,t)}{k}\fstop
$$
Using the evenness, we therefore obtain
$$
\Var(J(t))
=\frac{1}{\pi^2}\int_0^{k_0}\frac{1-B(k,t)}{k}\diff k+O(1)\fstop
$$
%All bounded errors below may depend on $k_0$, which remains fixed.

\smallskip
\noindent\emph{Step 2: control frequencies below $t^{-1}$.}
Take $t$ sufficiently large that $t^{-1}<k_0$. By
Lemma~\ref{lem:S-bounds} combined  with~\eqref{d:B}, we have
$$
0\le1-B(k,t)
=\frac{2\pi}{k}\big(S(k,0)-S(k,t)\big)
\le\pi tk\comma\qquad 0<k\le k_0\fstop
$$
It follows that
$$
0\le\int_0^{1/t}\frac{1-B(k,t)}{k}\diff k
\le\pi t\int_0^{1/t}\diff k=\pi\fstop
$$
Thus the  singularity at zero does not make the integral exploded.
%contributes only a bounded amount.

\smallskip
\noindent\emph{Step 3: extract the logarithm from the remaining band.}
On $[t^{-1},k_0]$, write
$$
\int_{1/t}^{k_0}\frac{1-B(k,t)}{k}\diff k
=\log(k_0t)-\int_{1/t}^{k_0}\frac{B(k,t)}{k}\diff k\fstop
$$
By Lemma~\ref{lem:B-bounds}, there is $c=c(k_0)>0$ such that
$B(k,t)\le e^{-ctk}$ on this interval. With $v=tk$, we obtain
$$
0\le\int_{1/t}^{k_0}\frac{B(k,t)}{k}\diff k
\le\int_1^{tk_0}\frac{e^{-cv}}{v}\diff v
\le\int_1^\infty\frac{e^{-cv}}{v}\diff v<\infty\fstop
$$
The last bound is independent of $t$. Combining the three steps gives
$\int_0^{k_0}(1-B(k,t))\diff k/k=\log t+O(1)$, and hence the result.
\end{proof}

\begin{rem}[Comparison with~\cite{Spo87}]\label{r:CS}
On top of the definition of the current~$J(t)$ already mentioned in~\eqref{d:RDJ}, our proof can be compared with the following points:
\begin{enumerate}

\item \emph{Infinite-particle limit.}
\cite[p.~670]{Spo87} first defines the expectations
for finite particle number and suppresses the subsequent
infinite-particle limit in the presentation.
We work directly with the stationary infinite system.
In particular, we avoid the formal subtraction
$\xi_t((0,\infty))-\xi_0((0,\infty))$,
where each of the two terms are almost surely infinite.

\item \emph{The variance formula.}
In the passage from~(2.7) to~(2.8),
\cite[pp.~672--673]{Spo87} uses test functions approaching
a delta function to obtain the current variance formula.
We provide the limiting argument needed to identify that formula
with the variance of the pathwise current~$J(t)$.
Lemma~\ref{l:WCV} identifies the $L^2$ limits of compactly
supported linear-statistic increments with the weighted currents.
Proposition~\ref{prop:direct-current-L2} then proves that the
compact half-line statistics converge to $J(t)$ in $L^2$,
giving the details in the passage to its variance.
\end{enumerate}
\end{rem}

\section{CLT for the current}\label{sec:CLT}

In this section, we prove the CLT for the current~$J(t)$ as $t \to \infty$. 

\begin{thm}[CLT for the current]\label{thm:CLT}
As $t\to\infty$, the following convergence in distribution holds:
$$
\frac{J(t)}{\sqrt{\Var(J(t))}}\ \Longrightarrow\ \mcN(0,1)\fstop
$$
Equivalently, by Theorem~\ref{thm:main},
$$
\frac{\pi\,J(t)}{\sqrt{\log t}}\ \Longrightarrow\ \mcN(0,1)\fstop
$$
\end{thm}

\subsection{Reduction to compactly supported two-time linear statistics}

Recall~\eqref{eq:sharp-current-cutoff}:
$$
f_R(x)=\chi(x/R)\1_{\{x>0\}}\comma\qquad
Z_R(t)=\int_{\R}f_R(x)
\bigl(\xi_t(\diff x)-\xi_0(\diff x)\bigr)\fstop
$$
Throught this section, the cut-off function $\chi$ is fixed. Constants below may depend
on this fixed cut-off function.

\begin{lem}\label{lem:CLT-reduction}
Set $R(t)\defi t^2$ and $X_t\defi Z_{R(t)}(t)$ for $t\ge2$.
Then
$$
\|X_t-J(t)\|_{L^2}\le\frac{\|\psi\|_{L^2}}{\sqrt t}
\longrightarrow0\comma\qquad
\frac{X_t-J(t)}{\sqrt{\Var(J(t))}}
\xrightarrow[t\to\infty]{L^2}0\fstop
$$
\end{lem}

\begin{proof}
The desired estimate follows from~\eqref{eq:sharp-current-error} in
Proposition~\ref{prop:direct-current-L2}. By substituting $R=t^2$, we obtain
 the first claim of the convergence. By Theorem~\ref{thm:main}, we obtain 
$\Var(J(t))\sim\pi^{-2}\log t$, which proves the second assertion.
\end{proof}
Using the preceding lemma, we can transfer the CLT from $X_t$ to $J(t)$: both variables~$X_t$ and $J(t)$ are centred, so the triangle
inequality provides
$$
\left|\sqrt{\Var(X_t)}-\sqrt{\Var(J(t))}\right|
\le\|X_t-J(t)\|_{L^2}\fstop
$$
By Slutsky's theorem,  the CLT can be transferred from $X_t$ to $J(t)$.
\subsection{Determinantal cumulants for a non-Hermitian kernel}\label{subsec:nonherm}

We first recall  several facts from the operator theory. %and then
For a compact operator $T$, its singular values $s_j(T)$ are
the eigenvalues of $(T^*T)^{1/2}$, counted with multiplicity.
The trace and Hilbert--Schmidt norms are
$$
 \|T\|_1\defi\sum_j s_j(T)\comma\qquad
 \|T\|_2\defi\left(\sum_j s_j(T)^2\right)^{1/2}\fstop
$$
The trace-class and
Hilbert--Schmidt operators are defined by the finiteness of these norms respectively.  
We write $\|T\|$ for
the usual operator norm. For trace-class $T$,
$\Tr T$ is the basis-independent sum of its diagonal matrix
entries in any orthonormal basis. We use the elementary bounds
\begin{equation}\label{eq:operator-ideal-facts}
 |\Tr T|\le\|T\|_1\comma\qquad
 \|H_1H_2\|_1\le\|H_1\|_2\|H_2\|_2\comma\qquad
 \|BTC\|_p\le\|B\|\,\|T\|_p\,\|C\|\quad(p=1,2)\comma
\end{equation}
where $H_1,H_2$ are Hilbert--Schmidt and $B,C$ are bounded.
The
Hilbert--Schmidt norm of an integral operator is the $L^2$
norm of its kernel. We also use cyclicity of the trace:
$\Tr(BT)=\Tr(TB)$ for bounded $B$ and trace-class $T$,
and $\Tr(H_1H_2)=\Tr(H_2H_1)$ for two Hilbert--Schmidt
operators. For functions, $\|\cdot\|_p$ denotes the usual
$L^p$ norm. 
Let $M_g$ denote the multiplication operator by a function~$g $. 
\begin{lem}[Integral Fredholm expansion]\label{lem:fredholm-series}
Let $\mathcal T$ be a finite nonempty set, and put
$E=\mathcal T\times\R$ with $\nu=\#\times\diff x$,
where $\#$ denotes counting measure on $\mathcal T$.
Let $A\subset E$ be compact. 
Suppose the chosen kernel $K$
is continuous in the spatial variables on each pair of
layers and defines a trace-class operator $K_A$ on
$L^2(A,\nu)$ as $K_A\defi \1_AK\1_A$, viz.,
\[
(K_Ah)(z):=\int_A K(z,w)h(w)\,\nu(dw),
\qquad z\in A.
\]
For every bounded measurable $g$ on $A$,
\begin{equation}\label{eq:fredholm-series-operator}
 \det(I+M_gK_A)
 =\sum_{n=0}^{\infty}\frac1{n!}
   \int_{A^n}\prod_{j=1}^n g(z_j)
   \det\big[K(z_i,z_j)\big]_{i,j=1}^n
   \diff\nu^{\otimes n}(z)\fstop
\end{equation}
The series is absolutely convergent.
\end{lem}

\begin{proof}
Put $\mcH=L^2(A,\nu)$, $V=\nu(A)$ and
$M=\sup_{A\times A}|K|$. 
Hadamard's inequality gives
$$
\frac1{n!}\int_{A^n}
\left|\prod_{i=1}^n g(z_i)\det[K(z_i,z_j)]_{i,j=1}^n\right|
\diff\nu^{\otimes n}(z)
\le\frac{n^{n/2}}{n!}(MV\|g\|_\infty)^n\fstop
$$
These bounds are summable, since $n!\ge(n/e)^n$. 

First suppose that $g$ is continuous and put $L=M_gK_A$.
By the  exterior-power expansion
\cite[Eq.~(3.3) and Theorem~3.2]{Simon1977}, we have
$$
\det(I+L)=\sum_{n=0}^{\infty}
\Tr\!\left(\bigwedge\nolimits^nL\right)\fstop
$$
Here $\bigwedge^nL$ acts on antisymmetric functions of $n$ variables.
We show how each of trace in this expansion becomes an integral in the RHS of~\eqref{eq:fredholm-series-operator}.
Fix $n\ge1$ and identify the tensor product $\mcH^{\otimes n}$ with
$\mcH_n\defi L^2(A^n,\nu^{\otimes n})$.
Let $\mcH_n^-$ be its closed subspace of antisymmetric functions:
these are the functions that change sign when any two arguments
are interchanged. The orthogonal projection onto this subspace is given by 
$$
(\mathcal A_nu)(z_1,\ldots,z_n)
=\frac1{n!}\sum_{\sigma\in\mathfrak S_n}
 \operatorname{sgn}(\sigma)
 u(z_{\sigma(1)},\ldots,z_{\sigma(n)})\comma
$$
where $\mathfrak S_n$ is the permutations group of
$\{1,\ldots,n\}$.

The operator $L^{\otimes n}$ applies $L$ in each coordinate, thus it commutes
with coordinate permutations as the same operator acts on each coordinate.  Hence it commutes with $\mathcal A_n$.
Its restriction to $\mcH_n^-$ is  $\bigwedge^nL$.
Consequently,
$$
T_n\defi\mathcal A_nL^{\otimes n}
       =L^{\otimes n}\mathcal A_n
$$
agrees with $\bigwedge^nL$ on $\mcH_n^-$ and is zero on
$(\mcH_n^-)^\perp$. It is trace class: the tensor product satisfies
$\|L^{\otimes n}\|_1=\|L\|_1^n<\infty$, and multiplication
by the bounded projection $\mathcal A_n$ preserves the trace class property.
We also have
$$
\Tr_{\mcH_n}(T_n)
=\Tr_{\mcH_n^-}\!\left(\bigwedge\nolimits^nL\right)\fstop
$$
Indeed, choose an orthonormal basis $(e_j)_j$ of $\mcH_n^-$
and an orthonormal basis $(v_k)_k$ of $(\mcH_n^-)^\perp$.
They form an orthonormal basis of $\mcH_n$, and $T_ne_j=(\bigwedge^nL)e_j$ whereas
$T_nv_k=0$. Since $T_n$ is trace class, 
$$
\begin{aligned}
 \Tr_{\mcH_n}(T_n)
 &=\sum_j\langle T_ne_j,e_j\rangle
   +\sum_k\langle T_nv_k,v_k\rangle
 &=\sum_j\left\langle
       \left(\bigwedge\nolimits^nL\right)e_j,e_j
     \right\rangle
 =\Tr_{\mcH_n^-}\!\left(\bigwedge\nolimits^nL\right)\fstop
\end{aligned}
$$

To compute the kernel of $T_n$, take
$\ell(z,w)\defi g(z)K(z,w)$, which is the kernel of $L$.
For $u\in\mcH_n$,
$$
(L^{\otimes n}u)(z)
=\int_{A^n}\prod_{i=1}^n\ell(z_i,w_i)\,
 u(w)\diff\nu^{\otimes n}(w)\fstop
$$
By applying $\mathcal A_n$ in the $z$ variables,  we thus obtain
$$
\begin{aligned}
H_n(z,w)
&:=\frac1{n!}\sum_{\sigma\in\mathfrak S_n}
  \operatorname{sgn}(\sigma)
  \prod_{i=1}^n\ell(z_{\sigma(i)},w_i)
&=\frac1{n!}\det[\ell(z_i,w_j)]_{i,j=1}^n
 =\frac1{n!}\det[g(z_i)K(z_i,w_j)]_{i,j=1}^n\fstop
\end{aligned}
$$

Since $g$ and $K$ are continuous, $H_n$ is continuous on
$A^n\times A^n$. Moreover, $A^n$ is a second-countable space
and $\nu^{\otimes n}$ restricted to it is a finite Borel measure.
The diagonal trace formula
\cite[Corollary~3.2]{Brislawn1991} therefore applies to $T_n$, which provides
$$
\begin{aligned}
\Tr\!\left(\bigwedge\nolimits^nL\right)
=\Tr(T_n)
 =\int_{A^n}H_n(z,z)\diff\nu^{\otimes n}(z) 
=\frac1{n!}\int_{A^n}
  \prod_{i=1}^n g(z_i)
  \det[K(z_i,z_j)]_{i,j=1}^n
  \diff\nu^{\otimes n}(z)\fstop
\end{aligned}
$$
In the last equality, we first set $w=z$ and then factor $g(z_i)$
out of the $i$th row of the determinant. This is the
$n$th term of~\eqref{eq:fredholm-series-operator}.
The $n=0$ term is $1$, so we have that for continuous $g$ 
\[
\begin{aligned}
\det(I+M_gK_A)
=1+\sum_{n=1}^{\infty}
\Tr\!\left(\bigwedge\nolimits^nL\right) 
=1+\sum_{n=1}^{\infty}\frac1{n!}
  \int_{A^n}
  \left(\prod_{i=1}^n g(z_i)\right)
  \det[K(z_i,z_j)]_{i,j=1}^n\,
  d\nu^{\otimes n}(z).
\end{aligned}
\]

For bounded measurable $g$, by Lusin's theorem and a continuous extension, 
we have a  continuous function $g_j$ on $A$ with $|g_j|\le\|g\|_\infty$ and
$g_j\to g$ almost everywhere. By the dominated convergence, we  have the strong
convergence of $M_{g_j}$ to $M_g$. Since $K_A$ is trace class,
$$
\|(M_{g_j}-M_g)K_A\|_1\longrightarrow0\fstop
$$
Indeed, approximate $K_A$ by a finite-rank operator $R$ and use
$$
\|(M_{g_j}-M_g)K_A\|_1
\le2\|g\|_\infty\|K_A-R\|_1
  +\|(M_{g_j}-M_g)R\|_1\fstop
$$
For fixed $R$, the second term tends to zero on its finitely many
range vectors; the first can be made arbitrarily small.
By the continuity of the determinant in trace norm
\cite[Theorem~3.5]{Simon1977}, we conclude the convergence of the operator
determinants. The dominated convergence again applies to each integral and
then to the series, by the common bound above.
This proves~\eqref{eq:fredholm-series-operator}.
\end{proof}

In the following proposition, we do not assume that the kernel is Hermitian. 
\begin{prop}\label{prop:fredholm}
Let $\mathcal T$ be a finite nonempty set, and put
$E=\mathcal T\times\R$ with $\nu=\#\times\diff x$.
Let $\Xi$ be a locally finite determinantal point process on $E$.
Assume that $(x,y)\mapsto K(s,x;r,y)$ is continuous on $\R^2$
for each $s,r\in\mathcal T$. For every compact
$A\subset\mathcal T\times\R$, assume that the operator on
$L^2(A,\nu)$ defined by
$$
 (K_Af)(z)\defi\int_A K(z,w)f(w)\,\nu(\diff w),
 \qquad z\in A,
$$
is trace class.
Let $F:E\to\mbbC$ be bounded and measurable with compact support and define 
$X_F\defi \sum_{z\in\Xi}F(z)$. Put $A\defi \supp F$ and
$K_A\defi \1_AK\1_A$. Then for every $\lambda\in\mbbC$,
$$
\mbbE\big[e^{\lambda X_F}\big]=\det_{L^2(A,\nu)}\!\big(I+M_{e^{\lambda F}-1}K_A\big)\comma
$$
where the determinant is the Fredholm determinant of the trace class operator
$M_{e^{\lambda F}-1}K_A$.
\end{prop}

\begin{proof}
This is the determinantal multiplicative functional identity, see, e.g., 
\cite[Eq.~(2.4)]{Borodin2011}. The statement there was given for discrete underlying spaces. We recall why it applies under our
continuous space setting. Put $G=e^{\lambda F}-1$. We then have
$$
e^{\lambda X_F}=\prod_{z\in\Xi\cap A}(1+G(z))\fstop
$$ 
Writing these subsets as
ordered tuples of distinct points introduces the factor $1/n!$.
Using the correlation densities, we can expand it as 
\begin{align*}
\mbbE e^{\lambda X_F}
&=\sum_{n=0}^{\infty}\frac1{n!}
  \int_{A^n}\prod_{j=1}^nG(z_j)\rho^{(n)}(z)
  \diff\nu^{\otimes n}(z)\\
&=\sum_{n=0}^{\infty}\frac1{n!}
  \int_{A^n}\prod_{j=1}^nG(z_j)
  \det[K(z_i,z_j)]_{i,j=1}^n\diff\nu^{\otimes n}(z)\\
&=\det_{L^2(A,\nu)}(I+M_GK_A)\fstop
\end{align*}
The absolute-convergence bound in Lemma~\ref{lem:fredholm-series},
applied to $|G|$, justifies interchange with expectation by Tonelli's
theorem. 
\end{proof}

We recall the definition of cumulants. Let $X$ be a real-valued
random variable. Its \emph{moment generating function} is defined as 
$$
 M_X(\lambda)\defi\mbbE[e^{\lambda X}]\comma
$$
for those real $\lambda$ for which the expectation is finite.
We assume that $X$ has a finite exponential moment:
$\mbbE[e^{\eta|X|}]<\infty$ for some $\eta>0$.
Under this assumption, the same formula
defines an analytic function of complex $\lambda$ for $|\lambda|<\eta$.
By taking the differentiation under the expectation, 
$$
 M_X^{(m)}(0)=\mbbE[X^m]\comma\qquad m\ge1\fstop
$$
Thus the derivatives of $M_X$ at zero is  the moments of $X$.

Since $M_X(0)=1$, the function $M_X$ has no zeros in a sufficiently
small disc around zero. On that disc,  we can define the
\emph{cumulant generating function}
$$
 \Lambda_X(\lambda)\defi\log M_X(\lambda)\comma
$$
using the unique analytic logarithm with $\Lambda_X(0)=0$.
The \emph{$m$th cumulant} of $X$ is defined as the derivative
$$
 \kappa_m(X)\defi\Lambda_X^{(m)}(0)
 =\left.\frac{\diff^m}{\diff\lambda^m}
       \log\mbbE[e^{\lambda X}]\right|_{\lambda=0}\comma
 \qquad m\ge1\fstop
$$
Equivalently, by the Taylor expansion,  for sufficiently small $|\lambda|$,
\begin{equation}\label{eq:cumulant-definition}
 \log\mbbE[e^{\lambda X}]
 =\Lambda_X(\lambda)
 =\sum_{m=1}^{\infty}\kappa_m(X)\frac{\lambda^m}{m!}\fstop
\end{equation}
In particular, the coefficient of $\lambda^m$ is
$\kappa_m(X)/m!$.

Using $M_X(0)=1$ and differentiating the logarithm, we have 
$$
 \begin{aligned}
 \kappa_1(X)
 &=M_X'(0)=\mbbE[X]\comma\\
 \kappa_2(X)
 &=M_X''(0)-\bigl(M_X'(0)\bigr)^2
   =\mbbE[X^2]-\bigl(\mbbE[X]\bigr)^2
   =\Var(X)\fstop
 \end{aligned}
$$
Finally, for real constants $a$ and $b$,
$M_{aX+b}(\lambda)=e^{b\lambda}M_X(a\lambda)$.
By taking logarithms near zero, we have 
$\Lambda_{aX+b}(\lambda)=b\lambda+\Lambda_X(a\lambda)$.
By taking the differentiating at zero, we have 
$$
 \kappa_1(aX+b)=a\kappa_1(X)+b\comma\qquad
 \kappa_m(aX+b)=a^m\kappa_m(X)\quad(m\ge2)\fstop
$$
Thus adding a constant affects only the first cumulant, while
multiplication by $a$ multiplies the $m$th cumulant by $a^m$.

\begin{prop} \label{prop:cumulant-trace}
Under the same setting as Proposition~\ref{prop:fredholm}, assume in addition that $F$ is real-valued. Then the cumulants $\kappa_m(X_F)$ satisfy that, for each $m\ge 1$,
\begin{equation}\label{eq:cumulant-trace}
\kappa_m(X_F)
=\sum_{\ell=1}^m(-1)^{\ell-1}\frac{m!}{\ell}\!\!\!\!\sum_{\substack{p_1+\cdots+p_\ell=m\\ p_j\ge 1}}
\frac{1}{p_1!\cdots p_\ell!}\,
\Tr_{L^2(A,\nu)}\Big(M_{F^{p_1}}K_A\,M_{F^{p_2}}K_A\cdots M_{F^{p_\ell}}K_A\Big)\comma
\end{equation}
where $A=\supp F$ and $K_A=\1_AK\1_A$.
\end{prop}

\begin{proof}
Put $A_\lambda=M_{e^{\lambda F}-1}K_A$. We have
$$
\|A_\lambda\|_1
\le\bigl(e^{|\lambda|\|F\|_\infty}-1\bigr)\|K_A\|_1
\longrightarrow0\qquad(\lambda\to0)\fstop
$$
By Plemelj's trace formula \cite[Eq.~(1.6)]{Simon1977} combined 
with Proposition~\ref{prop:fredholm}, we have that on a small disc around zero
$$
\log\mbbE e^{\lambda X_F}
=\log\det(I+A_\lambda)
=\sum_{\ell\ge1}\frac{(-1)^{\ell-1}}{\ell}
\Tr(A_\lambda^\ell)\fstop
$$ 
To identify its coefficients,
we expand each multiplier using
$e^{\lambda F}-1=\sum_{p\ge1}\lambda^pF^p/p!$. For fixed $\ell$,
$$
A_\lambda^\ell
=\sum_{p_1,\ldots,p_\ell\ge1}
\frac{\lambda^{p_1+\cdots+p_\ell}}{p_1!\cdots p_\ell!}
M_{F^{p_1}}K_A\cdots M_{F^{p_\ell}}K_A\fstop
$$
Each factor provides at least one power of $\lambda$. Thus the
coefficient of $\lambda^m$ uses only $\ell\le m$ and the positive
compositions $p_1+\cdots+p_\ell=m$. Comparing it  with
\eqref{eq:cumulant-definition}, 
we proved~\eqref{eq:cumulant-trace}.
\end{proof}

\subsection{Two-time extended sine kernel and a uniform cumulant bound}
\label{subsec:cumulant-bounds}

Fix $t>0$ and consider 
$$
 E_t\defi (\{0\}\times\R)\ \sqcup\ (\{t\}\times\R)
$$
with reference measure $\nu\defi\#\times\diff x$. We then define
$$
 \Xi_t\defi(\{0\}\times\xi_0)\cup(\{t\}\times\xi_t) \comma
$$
which has the correlation kernel~\eqref{eq:extended-sine-kernel},
by~\eqref{eq:extended-sine-correlations}.
Define
$$
 p(k)\defi\1_{\{|k|\le\pi\}}\comma
 \qquad q(k)\defi1-p(k)\fstop
$$
The Fourier multipliers are given by~\eqref{eq:kappa-KT}.
The multiplier $\kappa_{0,t}$ grows with $t$, which is inconvenient
for uniform operator estimates for later arguments. We remove this growth by setting
$\gamma_s\defi\exp(\pi^2s/2)$ for $s\ge0$. We then define {\it the gauge transformed kernel}
\begin{equation}\label{eq:gauge-kernel}
 \widetilde K_{\sin}(s,x;\,r,y)
 \defi\frac{\gamma_s}{\gamma_r}K_{\sin}(s,x;\,r,y)\comma
 \qquad s,r\ge0\fstop
\end{equation}
This gauge transformation preserves all correlation determinants
in~\eqref{eq:extended-sine-correlations}.
Indeed, for points $(s_i,x_i)$, it multiplies row $i$ by
$\gamma_{s_i}$ and column $i$ by $\gamma_{s_i}^{-1}$.
These factors cancel in the determinant. Thus
$\widetilde K_{\sin}$ defines the same determinantal field as
$K_{\sin}$ on any finite collection of time layers.

On the two layers $0$ and $t$, the Fourier multipliers of
$\widetilde K_{\sin}$ are
$\widetilde\kappa_{s,r}=(\gamma_s/\gamma_r)\kappa_{s,r}$.
Letting 
$$
 b_t(k)\defi\exp\!\left[-\frac t2(\pi^2-k^2)\right]\comma
 \qquad
 c_t(k)\defi\exp\!\left[-\frac t2(k^2-\pi^2)\right]\comma
$$
we obtain
\begin{equation}\label{eq:kappa-def}
 \begin{aligned}
 \widetilde\kappa_{0,0}(k)
   =\widetilde\kappa_{t,t}(k)&=p(k)\comma\\
 \widetilde\kappa_{0,t}(k)&=b_t(k)p(k)\comma\\
 \widetilde\kappa_{t,0}(k)&=-c_t(k)q(k)\fstop
 \end{aligned}
\end{equation}
Here $0\le b_t(k)\le1$ on $\{|k|\le\pi\}$ and
$0\le c_t(k)\le1$ on $\{|k|>\pi\}$. Consequently, all four
multipliers in~\eqref{eq:kappa-def} are bounded above by $1$, uniformly in $t>0$.

\subsubsection{Operator form in Fourier variables}
\label{subsec:operator-fourier}
We now turn the Fourier formula for the kernel into an operator
on the two time layers.

\smallskip
\noindent\emph{The Fourier projections.}
For $h\in L^2(\R)$, define $P$ and $Q$ by
\begin{equation} \notag %\label{eq:PQ-action}
 \widehat{Ph}(k)=p(k)\widehat h(k)\comma
 \qquad
 \widehat{Qh}(k)=q(k)\widehat h(k)\fstop
\end{equation}
Thus $P$ retains the frequencies $|k|\le\pi$ whereas
$Q=I-P$ retains their complement. Since their symbols~$p$ and $q$ are complementary real indicators,
Plancherel's theorem shows that $P$ and $Q$ are orthogonal
projections with $PQ=QP=0$ and $P+Q=I$.

For $h\in L^2(\R)$,
define
\begin{align} \label{e:BC}
 (B_th)(x)
 &\defi\frac1{2\pi}\int_{|k|\le\pi}
 e^{ikx}\exp\!\left[-\frac t2(\pi^2-k^2)\right]
                    \widehat h(k)\diff k\comma
%\notag % \label{eq:BC-action}
\\
 (C_th)(x)
 &\defi\frac1{2\pi}\int_{|k|>\pi}
 e^{ikx}\exp\!\left[-\frac t2(k^2-\pi^2)\right]
                    \widehat h(k)\diff k\fstop\notag
\end{align}
For each fixed $t>0$, both symbols belong to $L^2(\R)$.
Their products with $\widehat h$ are therefore in $L^1(\R)$,
so these integrals are also absolutely convergent. Their Fourier transforms are
$b_t p\widehat h$ and $c_t q\widehat h$, respectively.
The symbols have absolute value at most one. Thus, the Plancherel formula
gives the operator norm bounds below. Their supports also
show that $B_t$ acts within $PL^2$ and vanishes on $QL^2$;
$C_t$ acts within $QL^2$ and vanishes on $PL^2$. Hence
\begin{equation}\label{eq:BC-sectors}
 \|B_t\|\le1\comma\qquad \|C_t\|\le1\comma\qquad
 B_t=PB_tP\comma\qquad C_t=QC_tQ\fstop
\end{equation}

\smallskip
\noindent\emph{The two-time space and the current test.}
A function on $E_t$ is a pair $h=(h_0,h_t)$, where
$h_0(x)=h(0,x)$ and $h_t(x)=h(t,x)$. Thus
$$
 \mcH\defi L^2(E_t,\nu)
 \simeq L^2(\R)\oplus L^2(\R)\comma\qquad
 \|h\|_{\mcH}^2=\|h_0\|_2^2+\|h_t\|_2^2\fstop
$$
Define
\begin{equation}\notag %\label{eq:K-block}
 \mbfK_t\defi
 \begin{pmatrix}
 P & B_t\\
 -C_t & P
 \end{pmatrix}\fstop
\end{equation}
Since $\|P\|,\|B_t\|,\|C_t\|\le1$, the triangle inequality gives,
for $h=(h_0,h_t)\in\mcH$,
$$
\begin{aligned}
 \|\mbfK_th\|_{\mcH}
 &\le \|(Ph_0,Ph_t)\|_{\mcH}
      +\|(B_th_t,-C_th_0)\|_{\mcH}
 &\le 2\|h\|_{\mcH}\fstop
\end{aligned}
$$
Taking the supremum over $\|h\|_{\mcH}=1$ yields
$\|\mbfK_t\|\le2$.
Its blocks have the symbols in~\eqref{eq:kappa-def}, so its
kernel is~the gauged kernel~$\widetilde K_{\sin}$. Note that this matrix is not 
self-adjoint in general. 

Define the test function $F_R:E_t\to\R$ as
\begin{equation}\notag %\label{eq:F-def}
 F_R(0,x)\defi-f_R(x)\comma
 \qquad
 F_R(t,x)\defi f_R(x)\fstop
\end{equation}
The two signs record the difference between the configurations
at times $t$ and $0$:
$$
 \sum_{z\in\Xi_t}F_R(z)
 =-\int_{\R}f_R(x)\xi_0(\diff x)
   +\int_{\R}f_R(x)\xi_t(\diff x)
 =Z_R(t)\fstop
$$
Writing $M_g$ for multiplication by $g$ on 
$L^2$ space, multiplication by $F_R^r$ acts separately
on each layer. Its matrix is
\begin{equation} \notag %\label{eq:MF-block}
 M_{F_R^r}
 =
 \begin{pmatrix}
 (-1)^r M_{f_R^r} & 0\\
 0 & M_{f_R^r}
 \end{pmatrix}\comma\qquad
 \|M_{F_R^r}\|=\|f_R\|_\infty^r\le1\fstop
\end{equation}
Here the last inequality uses $0\le f_R\le1$.

\begin{lem}
\label{lem:local-trace-class}
For a compact $A\subset E_t$, let $\Pi=M_{\1_A}$ on
$\mcH$, and let $K_A$ be the restriction of
$\Pi\mbfK_t\Pi$ to $L^2(A,\nu)$.
Then $K_A$ is trace class and its chosen kernel is continuous
in the spatial variables on each pair of time layers.
The same assertions hold for any finite
collection of distinct time layers of the gauged extended sine
kernel.
\end{lem}

\begin{proof}
The Fourier multipliers of~\(P\), \(B_t\), and \(-C_t\) are \(p\), \(b_tp\), and \(-c_tq\), respectively. Their pointwise absolute values are therefore \(p\), \(b_tp\), and \(c_tq\).
They are nonnegative, bounded and integrable. Thus their
inverse Fourier transforms are continuous.

Let $(D, d) \in \{(P,p), (B_t, b_tp), (C_t, c_tq)\}$
and let $I$ be a compact
interval. The square-root $\sqrt d$ belongs to~$L^2$.
If $L$ is its inverse Fourier transform, the Plancherel formula provides
$$
\|D^{1/2}M_{\1_I}\|_2^2
=\int_I\int_{\R}|L(x-y)|^2\diff x\diff y
=\frac{|I|}{2\pi}\int_{\R}d(k)\diff k<\infty\fstop
$$
Hence
$M_{\1_I}DM_{\1_I}=(D^{1/2}M_{\1_I})^*(D^{1/2}M_{\1_I})$
is trace class, by~\eqref{eq:operator-ideal-facts}.

Choose a compact interval $I$ containing the spatial part of $A$
on every layer. We define $M_I\defi M_{\1_I}$ and 
$$
 P_I\defi M_IPM_I\comma\qquad
 B_{t,I}\defi M_IB_tM_I\comma\qquad
 C_{t,I}\defi M_IC_tM_I\fstop
$$
The above factorisation shows that all three operators are of 
trace class on $L^2(\R)$. We now show that the operator on
$\mcH=L^2(\R)\oplus L^2(\R)$ given by
$$
 T_I\defi
 \begin{pmatrix}
 P_I & B_{t,I}\\
 -C_{t,I} & P_I
 \end{pmatrix}
$$
is also  trace class.

Let $U_0h\defi(h,0)$ and $U_th\defi(0,h)$ be the inclusions
of $L^2(\R)$ into the two components of $\mcH$.
Their adjoints are the coordinate projections:
$U_0^*(h_0,h_t)=h_0$ and $U_t^*(h_0,h_t)=h_t$.
Thus
$$
 T_I
 =U_0P_IU_0^*+U_0B_{t,I}U_t^*
  -U_tC_{t,I}U_0^*+U_tP_IU_t^*\fstop
$$
Each summand is  trace class since its middle factor is trace
class and its two outer (left and right) factors are bounded. All the inclusions
and projections have norm one, so~\eqref{eq:operator-ideal-facts}
and the triangle inequality give
$$
 \|T_I\|_1
 \le 2\|P_I\|_1+\|B_{t,I}\|_1+\|C_{t,I}\|_1<\infty\fstop
$$

Finally, let $J_A:L^2(A,\nu)\to\mcH$ be extension by zero;
its adjoint $J_A^*$ restricts a function to $A$.
Because $I$ contains the spatial part of $A$ on both layers,
the extra cut-offs in $T_I$ have no effect after this restriction.
Consequently,
$$
 K_A=J_A^*T_IJ_A\comma\qquad
 \|K_A\|_1\le\|T_I\|_1<\infty\fstop
$$
This proves that $K_A$ is a trace class operator on $L^2(A,\nu)$, which completes the case of two time layers. 

The extension to the multi-time layer case is  straightforward: 
for time-points $t_0<\cdots<t_r$, choose again a compact interval $I$ containing the spatial part of $A$ on every layer.
The operator restricted spatially to $I$ on each layer has
diagonal blocks $P_I$ and off-diagonal blocks $B_{h,I}$ or
$-C_{h,I}$, where $h=|t_a-t_b|>0$.
The  same argument of the factorisation in the case of two time layer applies to each positive time gap, and   every block is trace class. Since there are finitely many
blocks, the whole operator is trace class. The continuity of the kernel on each pair of time
layers follows from the same Fourier inversion argument.
This proves the final assertion and completes the proof.
\end{proof}

\smallskip
Define
$$
 \Lambda_R\defi\supp f_R\comma\qquad
 A_R\defi\supp F_R
 =(\{0\}\times\Lambda_R)
     \sqcup(\{t\}\times\Lambda_R)\fstop
$$
The projection $\Pi=M_{\1_A}$ for a compact $A\subset E_t$
sets a function equal to zero outside $A$. 

\begin{lem}[Localisation of traces]\label{lem:trace-localisation}
Fix $t>0$, let $A\subset E_t$ be compact, and put
$\Pi=M_{\1_A}$ on $\mcH$. We identify $L^2(A,\nu)$
with the subspace of $\mcH$ consisting of functions that
vanish outside $A$, by extension by zero. Let
$$
 K_A=\left.\Pi\mbfK_t\Pi\right|_{L^2(A,\nu)}\fstop
$$
Let $\ell\ge1$, let $g_1,\ldots,g_\ell$ be bounded
measurable functions supported in $A$, and set
$$
 T\defi M_{g_1}\mbfK_t\cdots M_{g_\ell}\mbfK_t\fstop
$$
Then $T\Pi$ is trace class and
\begin{equation}\label{eq:trace-localisation}
 \Tr_{L^2(A,\nu)}
 \big(M_{g_1}K_A\cdots M_{g_\ell}K_A\big)
 =\Tr_{\mcH}(T\Pi)
\ 
. \end{equation}
If $\ell\ge2$, then $T$
 is trace class and
\begin{equation}\label{eq:trace-full-space}
 \Tr_{L^2(A,\nu)}
 \big(M_{g_1}K_A\cdots M_{g_\ell}K_A\big)
 =\Tr_{\mcH}(T)\fstop
\end{equation}
For \(\ell=1\), the operator \(M_{g_1}\mathbf K_t\Pi\) is also trace class, and
\begin{align} \label{e:L1}
\operatorname{Tr}_{L^2(A,\nu)}(M_{g_1}K_A)
=
\operatorname{Tr}_{\mathcal H}
\bigl(M_{g_1}\mathbf K_t\Pi\bigr).
\end{align}
The same conclusions hold for  any
finite collection of  time layers.
\end{lem}

\begin{proof}
Write $D_j=M_{g_j}$. As $g_j$ is supported in $A$, we have
$\Pi D_j=D_j\Pi=D_j$. Thus the extension by zero of the restricted
product can be writte as 
$$
D_1\Pi\mbfK_t\Pi\cdots D_\ell\Pi\mbfK_t\Pi=T\Pi\fstop
$$
All the projection operators are absorbed by adjacent multiplication factors, and 
this extension has the same
trace as the restricted operator, proving
\eqref{eq:trace-localisation}. We also have $\Pi T=T$.

We next show that each operator $D_j\mbfK_t$ is a Hilbert--Schmidt operator.
For $s,r\in\{0,t\}$, write the corresponding kernel block as
$$
 \widetilde K_{\sin}(s,x;r,y)=L_{s,r}(x-y)\fstop
$$ 
Its Fourier transform is 
$\widetilde\kappa_{s,r}$ in~\eqref{eq:kappa-def}.
These multipliers are bounded and integrable for fixed $t>0$,
in particular belong to $L^2(\R)$. By the Plancherel  theorem, we therefore have 
$$
 \|L_{s,r}\|_{L^2}^2
 =\frac1{2\pi}\int_{\R}|\widetilde\kappa_{s,r}(k)|^2\diff k
 <\infty\fstop
$$

The multiplication by $D_j=M_{g_j}$ acts on the output variable.
Thus the $(s,r)$ block of $D_j\mbfK_t$ has the integral kernel
$g_j(s,x)L_{s,r}(x-y)$.
For a fixed $x$, by the change of variable $u=x-y$, we have 
$$
 \int_{\R}|L_{s,r}(x-y)|^2\diff y
 =\int_{\R}|L_{s,r}(u)|^2\diff u
 =\|L_{s,r}\|_{L^2}^2\fstop
$$
By using this identity and Tonelli's theorem, we obtain
$$
\begin{aligned}
 &\int_{\R^2}|g_j(s,x)L_{s,r}(x-y)|^2\diff x\diff y
 &
 =\int_{\R}|g_j(s,x)|^2
       \left(\int_{\R}|L_{s,r}(x-y)|^2\diff y\right)\diff x
 &
 =\|g_j(s,\cdot)\|_{L^2}^2\|L_{s,r}\|_{L^2}^2<\infty\fstop
\end{aligned}
$$
The last quantity is finite because $g_j$ is bounded and
supported in the compact set $A$. This proves that every
block is a Hilbert--Schmidt operator. On the two-time space, the squared
Hilbert--Schmidt norm is the sum over the four blocks:
$$
 \|D_j\mbfK_t\|_2^2
 =\sum_{s,r\in\{0,t\}}
   \|g_j(s,\cdot)\|_{L^2}^2\|L_{s,r}\|_{L^2}^2<\infty\fstop
$$

Suppose now that $\ell\ge2$. We group the first two factors in
$T=(D_1\mbfK_t)\cdots(D_\ell\mbfK_t)$.
Their product is trace class since both are Hilbert--Schmidt.
The multiplication by the remaining bounded factors preserves the trace
class property. More explicitly,~\eqref{eq:operator-ideal-facts} gives
$$
 \|T\|_1
 \le \|D_1\mbfK_t\|_2\,\|D_2\mbfK_t\|_2
       \prod_{j=3}^{\ell}\|D_j\mbfK_t\|<\infty\comma
$$
where the empty product is $1$ when $\ell=2$.

We can now apply the trace cyclicity since $T$ is trace class and
$\Pi$ is bounded. Recall also that the leftmost multiplier
$D_1$ makes $T$ vanish outside $A$:
$$
 \Pi T
 =(\Pi D_1)\mbfK_tD_2\mbfK_t\cdots D_\ell\mbfK_t
 =D_1\mbfK_tD_2\mbfK_t\cdots D_\ell\mbfK_t
 =T\fstop
$$
Consequently,
$$
 \Tr(T\Pi)=\Tr(\Pi T)=\Tr(T)\fstop
$$
The first equality is the cyclicity of trace, and the second uses $\Pi T=T$.
Combining this with~\eqref{eq:trace-localisation} proves
\eqref{eq:trace-full-space}.

For $\ell=1$, the estimate above only establishes that $T$ is
Hilbert--Schmidt. The first paragraph already supplies the
trace-class operator $T\Pi$, so we obtain~\eqref{e:L1}. 

For any finite collection of time layers, the same
argument applies with a finite sum over all pairs of layers, so we omit the proof. 
\end{proof}

\subsubsection{Fourier and commutator estimates}

We will control higher cumulants.

\begin{lem}\label{lem:fourier-powers}
For each integer $m\ge1$, there is $C_m<\infty$, depending only on
$m$ and the fixed cut-off $\chi$, such that for every $R\ge1$,
\begin{equation}\label{eq:hat-power-bound}
\big|\widehat{f_R^m}(k)\big|
\le
\begin{cases}
C_m\min\{R,|k|^{-1}\},&k\ne0,\\
2R,&k=0.
\end{cases}
\end{equation}
\end{lem}

\begin{proof}
Since $f_R^m=\chi_R^m\1_{\{x>0\}}$ is supported in $[0,2R]$
and takes values in $[0,1]$, its $L^1$ norm is at most $2R$.
This proves $|\widehat{f_R^m}(k)|\le2R$ for every $k$.

For the second estimate, fix $k\ne0$ and take integration by parts on
$[0,2R]$, on which $\chi_R^m$ is smooth. Its endpoint values are
$\chi_R(0)^m=1$ and $\chi_R(2R)^m=0$, so
$$
ik\widehat{f_R^m}(k)
=1+\int_0^{2R}e^{-ikx}(\chi_R^m)'(x)\diff x\fstop
$$
The boundary term $1$ records the jump of $f_R^m$ at the origin.
By the chain rule and the change of variables $x=Ry$,
$$
\int_0^{2R}|(\chi_R^m)'(x)|\diff x
=m\int_0^2\chi(y)^{m-1}|\chi'(y)|\diff y
\le m\|\chi'\|_{L^1}\fstop
$$
It follows that
$$
|\widehat{f_R^m}(k)|
\le\frac{1+m\|\chi'\|_{L^1}}{|k|}\fstop
$$
Combining this with the $L^1$ estimate, we obtain 
\eqref{eq:hat-power-bound}.
\end{proof}

For linear maps \(A,B:V\to V\) on a vector space~$V$, their commutator is defined as 
\[
[A,B]\coloneqq AB-BA.
\]
\begin{lem}
\label{lem:comm-HS}\label{lem:offdiag-HS}
Let $g\in L^\infty(\R)$ with $\widehat g\in L^2(\R)$.
Then $PM_gQ$, $QM_gP$ and $[P,M_g]$ are Hilbert--Schmidt, and
\begin{equation}\label{eq:PMQ-HS}
 \|PM_gQ\|_2^2=\|QM_gP\|_2^2
 =\frac1{4\pi^2}\int_{\R}\min(|k|,2\pi)
                      |\widehat g(k)|^2\diff k\fstop
\end{equation}
Consequently, we have the commutator estimate
\begin{equation}\label{eq:comm-HS}
 \|[P,M_g]\|_2^2
 =\frac1{2\pi^2}\int_{\R}\min(|k|,2\pi)
                      |\widehat g(k)|^2\diff k\comma
\end{equation}
and, in particular,
\begin{equation}\label{eq:PMQ-comm}
 \|PM_gQ\|_2\le\|[P,M_g]\|_2\comma\qquad
 \|QM_gP\|_2\le\|[P,M_g]\|_2\fstop
\end{equation}
\end{lem}

\begin{proof}
Under the map
$h\mapsto(2\pi)^{-1/2}\widehat h$, the multiplication by $g$
has kernel $(2\pi)^{-1}\widehat g(k-k')$.
For fixed $g$, this formula is first verified on Schwartz functions and then can be extended to elements in $L^2(\R)$ by 
using the density of $\mcS(\R)$ in $L^2(\R)$ and the bound
$$
 \|M_gh\|_{L^2}\le\|g\|_\infty\|h\|_{L^2}\fstop
$$
 Since $P$ and $Q$ multiply
by $p$ and $q$ in Fourier variables, the kernel of $PM_gQ$ is given as 
$$
 (k,k')\longmapsto
 \frac1{2\pi}p(k)q(k')\widehat g(k-k')\fstop
$$
Its squared Hilbert--Schmidt norm is the integral of the
squared absolute value of this kernel. With $r=k-k'$, it is
$$
 \frac1{4\pi^2}\int_{\R}|\widehat g(r)|^2
       \left(\int_{\R}p(k)q(k-r)\diff k\right)\diff r\fstop
$$
The inner integral is the length of the part of
$[-\pi,\pi]$ lying outside $[r-\pi,r+\pi]$.
For $|r|\le2\pi$, this length is $|r|$, and for
$|r|\ge2\pi$, the intervals have disjoint interiors and
the length is $2\pi$. Thus it is equal to  $\min(|r|,2\pi)$.
The same calculation for $QM_gP$ gives the length of
$[r-\pi,r+\pi]\setminus[-\pi,\pi]$, which is equal
to the first length. This proves~\eqref{eq:PMQ-HS}.
The integrals are finite because $\widehat g\in L^2$
and $\min(|r|,2\pi)\le2\pi$.

Finally, on the space~$PL^2(\R)\oplus QL^2(\R)$,
$$
 [P,M_g]=
 \begin{pmatrix}0&PM_gQ\\-QM_gP&0\end{pmatrix}\fstop
$$
The squared Hilbert--Schmidt norm of this block matrix is
the sum of the squared norms of its two nonzero blocks.
Hence~\eqref{eq:comm-HS} and~\eqref{eq:PMQ-comm} follow
from~\eqref{eq:PMQ-HS}.
\end{proof}

\begin{lem}[A Schatten trace inequality]\label{lem:schatten}
Let $B_1,\dots,B_r$ be bounded operators on a Hilbert space and let $H_1,H_2$ be Hilbert--Schmidt operators.
Then $H_1B_1H_2B_2\cdots B_r$ is trace class and
\begin{equation}\notag %\label{eq:schatten}
\big|\Tr(H_1B_1H_2B_2\cdots B_r)\big|
\le \|H_1\|_2\,\|H_2\|_2\,\prod_{j=1}^r \|B_j\|\fstop
\end{equation}
\end{lem}

\begin{proof}
We group the product into two factors, each containing one of the
Hilbert--Schmidt operators. Then, the bounds in
\eqref{eq:operator-ideal-facts} make both groups Hilbert--Schmidt
with norms bounded by the corresponding Hilbert--Schmidt norm
times the norms of their bounded factors. Their product is trace
class, and the same bounds give
$$
\big|\Tr(H_1B_1H_2B_2\cdots B_r)\big|
\le\|H_1B_1\|_2\|H_2B_2\cdots B_r\|_2
\le\|H_1\|_2\|H_2\|_2\prod_j\|B_j\|\fstop
$$
Here note that an empty product of bounded factors is interpreted as the identity.
\end{proof}

\begin{lem}\label{lem:comm-logR}
There exists $C<\infty$ such that for all $R\ge 2$,
\begin{equation}\label{eq:comm-logR}
\|[P,M_{f_R}]\|_2^2\le C\,\log R\fstop
\end{equation}
More generally, if $g$ is bounded, $\|g\|_\infty\le1$,
and $[P,M_g]$ is Hilbert--Schmidt, then for every integer $m\ge1$,
\begin{equation}\label{eq:comm-power}
 \|[P,M_{g^m}]\|_2\le m\,\|[P,M_g]\|_2\fstop
\end{equation}
\end{lem}

\begin{proof}
Insert the Fourier bounds of Lemma~\ref{lem:fourier-powers}
with $m=1$ into~\eqref{eq:comm-HS}. Separating the domain of the integral into
$0<|k|\le R^{-1}$, $R^{-1}<|k|\le1$ and $|k|>1$, we have
$$
 \|[P,M_{f_R}]\|_2^2
 \le C\left(
 R^2\int_0^{1/R}k\diff k
 +\int_{1/R}^1\frac{\diff k}{k}
 +\int_1^\infty\frac{\diff k}{k^2}\right)
 \le C(1+\log R)\fstop
$$
For $R\ge2$ this proves~\eqref{eq:comm-logR}.

For~\eqref{eq:comm-power}, the identity
$[P,M_{uv}]=M_u[P,M_v]+[P,M_u]M_v$ gives
$$
 [P,M_{g^m}]
 =\sum_{r=0}^{m-1}M_{g^r}[P,M_g]M_{g^{m-1-r}}\fstop
$$
The operator norm of every multiplication factor in this sum is bounded above by $1$. Therefore, the Hilbert--Schmidt
norm of each summand is bounded by $\|[P,M_g]\|_2$. 
Thus, we have
$$
\begin{aligned}
 \|[P,M_{g^m}]\|_2
 &\le \sum_{r=0}^{m-1}
       \|M_{g^r}[P,M_g]M_{g^{m-1-r}}\|_2
 &\le \sum_{r=0}^{m-1}
       \|M_{g^r}\|\,\|[P,M_g]\|_2\,\|M_{g^{m-1-r}}\|
 &\le m\,\|[P,M_g]\|_2\fstop
\end{aligned}
$$
This proves~\eqref{eq:comm-power}.
\end{proof}

\subsubsection{Mixed-time traces on finitely many time layers}
\label{subsec:mixed-time-traces}
We formulate the trace estimate for a finite collection of 
time points.
Fix an integer $r\ge1$ and times $0=t_0<t_1<\cdots<t_r$.
For $0\le a,b\le r$, define 
$$
 (K_{t_a,t_b}u)(x)\defi
 \int_\R\widetilde K_{\sin}(t_a,x;\,t_b,y)u(y)\diff y\comma
 \qquad u\in L^2(\R)\fstop
$$
Here $\widetilde K_{\sin}$ is the gauged kernel
introduced in~\eqref{eq:gauge-kernel}; the row time $t_a$ is the output
layer and the column time $t_b$ is the input layer.
By Equations~\eqref{eq:ext-sine-FT} and~\eqref{eq:gauge-kernel}, we have 
$$
 \widehat{K_{t_a,t_b}u}(k)
 =\frac{\gamma_{t_a}}{\gamma_{t_b}}
       \kappa_{t_a,t_b}(k)\widehat u(k)\fstop
$$
By substituting~\eqref{eq:kappa-KT},  these operators are identified as
\begin{equation}\label{eq:multitime-gauge-blocks}
 K_{t_a,t_b}=
 \begin{cases}
 P,&a=b,\\
 B_{t_b-t_a},&a<b,\\
 -C_{t_a-t_b},&a>b.
 \end{cases}
\end{equation}
For a positive time gap $h$, the operators $B_h$ and $C_h$ are
those defined in Section~\ref{subsec:operator-fourier}, with $t$
replaced by $h$. Their Fourier multipliers are $b_h p$ and $c_h q$,
respectively. In particular, \eqref{eq:BC-sectors} gives
$$
 B_h=PB_hP\comma\qquad C_h=QC_hQ\comma\qquad
 \|P\|,\|B_h\|,\|C_h\|\le1\fstop
$$
Thus a block with $a\le b$ acts within $PL^2(\R)$ whereas
a block with $a>b$ acts within $QL^2(\R)$.

A \emph{cyclic time word of length $\ell$} is a tuple of layer
indices
$$
 \mathbf a=(a_1,\ldots,a_\ell)\in\{0,\ldots,r\}^{\ell}\comma
 \qquad a_{\ell+1}\defi a_1\fstop
$$
The convention $a_{\ell+1}=a_1$ closes the sequence of kernel
blocks by returning to its starting layer. The word is
\emph{nonconstant} if it visits at least two different layers.

\begin{lem}[Mixed-time trace bound]\label{lem:two-interfaces}
Let $R\ge2$, let $\mathbf a$ be a nonconstant cyclic time word,
and let $\mathbf p=(p_1,\ldots,p_\ell)$ consist of positive integers.
In particular, $\ell\ge2$. Set
$$
 L_j\defi K_{t_{a_j},t_{a_{j+1}}}\comma\qquad
 T_{\mathbf a,\mathbf p}\defi
 M_{f_R^{p_1}}L_1\cdots M_{f_R^{p_\ell}}L_\ell\fstop
$$
Then $T_{\mathbf a,\mathbf p}$ is trace class on $L^2(\R)$,
and its trace agrees with the corresponding localised trace
in Lemma~\ref{lem:trace-localisation}. Moreover,
\begin{equation}\label{eq:two-interfaces-bound}
 \big|\Tr_{L^2(\R)}(T_{\mathbf a,\mathbf p})\big|
 \le
 \left(\max_{1\le j\le\ell}
          \|[P,M_{f_R^{p_j}}]\|_2\right)^2\fstop
\end{equation}
This estimate is uniform in the number of  time-points,
their values, and their  gaps.
\end{lem}

\begin{proof}
Write $g_j\defi f_R^{p_j}$, so that
$T_{\mathbf a,\mathbf p}=M_{g_1}L_1\cdots M_{g_\ell}L_\ell$.
Each operator $M_{g_j}L_j$ is a Hilbert--Schmidt operator due to the blockwise estimate
in the proof of Lemma~\ref{lem:trace-localisation}, and the product is trace class.
The identification to its localised trace can be shown by applying the latter assertion of Lemma~\ref{lem:trace-localisation}  regarding finite-layer 
to 
$$
 h_j(s,x)\defi\1_{\{s=t_{a_j}\}}g_j(x)\fstop
$$
Then the resulting product on $\bigoplus_{a=0}^r L^2(\R)$ has
$T_{\mathbf a,\mathbf p}$ as its only possibly non-zero diagonal
block. Its full space trace is therefore
$\Tr(T_{\mathbf a,\mathbf p})$, and the same lemma gives the claimed
agreement with the localised trace.
In the following argument we prove the remainder of the statement. 

\medskip
\noindent\emph{Step 1: insert the frequency projections.}
For each kernel block, we define
$$
 S_j\defi
 \begin{cases}
 P,&a_j\le a_{j+1},\\
 Q,&a_j>a_{j+1},
 \end{cases}
 \qquad S_0\defi S_\ell\fstop
$$
By~\eqref{eq:multitime-gauge-blocks} and~\eqref{eq:BC-sectors},
$L_j=S_jL_jS_j$. Thus $L_j$ vanishes on the complementary
subspace and maps $S_jL^2(\R)$ into itself. Define
$$
 H_j\defi S_{j-1}M_{g_j}S_j\fstop
$$
For $2\le j\le\ell$, the projections~$S_j$ before and after $M_{g_j}$ can
be inserted because
$$
 L_{j-1}M_{g_j}L_j
 =L_{j-1}(S_{j-1}M_{g_j}S_j)L_j
 =L_{j-1}H_jL_j\fstop
$$
At the closing position, $L_\ell=L_\ell S_\ell$ gives
$T_{\mathbf a,\mathbf p}=T_{\mathbf a,\mathbf p}S_\ell$.
Since this operator is trace class, the cyclicity of trace implies
$$
\begin{aligned}
 \Tr(T_{\mathbf a,\mathbf p})
 &=\Tr(T_{\mathbf a,\mathbf p}S_\ell)
   =\Tr(S_\ell T_{\mathbf a,\mathbf p})
 &=\Tr(H_1L_1\cdots H_\ell L_\ell)\fstop
\end{aligned}
$$
The initial projection $S_\ell=S_0$ completes $H_1$,
and the remaining projections are inserted as above.

\medskip
\noindent\emph{Step 2.}
Since the word is cyclic, we have the cancellation:
$$
 \sum_{j=1}^{\ell}(a_{j+1}-a_j)=0\fstop
$$
Since the word is nonconstant, not all summand vanish, which means that 
there is at least one increment which is  positive and at least one increment
is negative. The case~$a_j<a_{j+1}$ gives a block $B_h$, hence $S_j=P$; the case~$a_j >a_{j+1}$ gives a block $-C_{h'}$, 
hence $S_j=Q$. Consequently, the cyclic sequence
$S_1,\ldots,S_\ell$ contains both $P$ and $Q$.
Starting at a $P$ entry and going once around the sequence,
we must enter the $Q$ entries and then return to $P$.
There are therefore two distinct positions $i<j$ with
$$
 S_{i-1}\ne S_i\comma\qquad S_{j-1}\ne S_j\fstop
$$
At either position $k\in\{i,j\}$,
$$
 H_k\in\{PM_{g_k}Q,\ QM_{g_k}P\}\comma\qquad
 \|H_k\|_2\le\|[P,M_{g_k}]\|_2\fstop
$$
The last inequality is~\eqref{eq:PMQ-comm}. 
For every $k$, we also have
$$
 \|H_k\|\le\|g_k\|_\infty\le1\comma\qquad
 \|L_k\|\le1\fstop
$$

\medskip
\noindent\emph{Step 3:  the trace estimate.}
Collect the factors before $H_i$, between $H_i$ and $H_j$,
and after $H_j$  into products $U,V,W$, respectively, with 
their original order kept. Then
$$
 H_1L_1\cdots H_\ell L_\ell=UH_iVH_jW\comma\qquad
 \|U\|,\|V\|,\|W\|\le1\fstop
$$
 The two groups $UH_iV$
and $H_jW$ are Hilbert--Schmidt. Thus, by using the Hilbert–Schmidt norm estimate and
Lemma~\ref{lem:schatten}, we have 
$$
\begin{aligned}
 |\Tr(T_{\mathbf a,\mathbf p})|
 &\le\|UH_iV\|_2\,\|H_jW\|_2
 &\le\|H_i\|_2\,\|H_j\|_2
 &\le\left(\max_{1\le k\le\ell}
                \|[P,M_{f_R^{p_k}}]\|_2\right)^2\fstop
\end{aligned}
$$
This proves~\eqref{eq:two-interfaces-bound}. The right-hand side of~\eqref{eq:two-interfaces-bound}
depends only on the test functions. Hence the estimate holds  independently of the number of observation
times and the sizes of their gaps.
The proof is complete.
\end{proof}

\subsubsection{One-time cumulant bound}
\label{subsec:sine-cumulants}

We next control the constant-time terms in the cumulant
expansion. At a single time, since the kernel is the orthogonal
projection $P$, we can apply the general determinant estimate of
Breuer and Duits,which gives a bound of  the centred logarithm
of the moment generating function by $\|[P,M_f]\|_2^2$.
We therefore need to check their hypothesis of the trace-class property.  We then extract
the cumulants by Cauchy's coefficient estimate.

Let the random variable $\xi$ have its law sine$_2$  under $\mathbb P$ and, for a bounded
compactly supported function~$f$, define
$$
 \mcN(f)\defi\int_{\R}f(x)\,(\xi(\diff x)-\diff x)\fstop
$$
This is the one-time version of $\mcN_t(f)$ defined earlier.

\begin{prop}\label{prop:sine-cumulant-bound}
Let $f$ be real-valued, bounded and compactly supported, with $0\le f\le 1$ and $[P,M_f]$ Hilbert--Schmidt. Then, for every integer $m\ge 3$,
\begin{equation}\label{eq:sine-cumulant-bound}
|\kappa_m(\mcN(f))|\le m!\,C^m\,\|[P,M_f]\|_2^2\comma
\end{equation}
for a universal constant $C<\infty$.
In particular, for $f=f_R$ and $R\ge2$,
$$
|\kappa_m(\mcN(f_R))|\le m!\,C^m\,\log R\fstop
$$
\end{prop}

\begin{proof}
We apply Breuer--Duits \cite[Lemma~3.2]{BreuerDuits2014} with
$K=P$ and $h=M_f$. 
The cited lemma needs that $M_f$ is bounded and both $M_fP$
and $PM_f$ are trace class. The boundedness follows from
$\|M_f\|=\|f\|_\infty\le1$. We now verify the trace-class
property of these two products.

Choose $L>0$ with $\supp f\subset[-L,L]$, and let
$Uu=(2\pi)^{-1/2}\widehat u$. The kernel of $M_fPU^{-1}$ expands as follows:
$$
\frac{f(x)e^{ikx}\1_{[-\pi,\pi]}(k)}{\sqrt{2\pi}}
=\sum_{n=0}^{\infty}
\frac{i^nx^nf(x)k^n\1_{[-\pi,\pi]}(k)}{\sqrt{2\pi}\,n!}\fstop
$$
Each term is rank one, and its trace norm is bounded as 
\[
\left\|\frac{i^nx^nf(x)}{\sqrt{2\pi}\,n!}\right\|_{L^2(dx)}
\left\|k^n\mathbf1_{[-\pi,\pi]}(k)\right\|_{L^2(dk)} \le \frac{\|f\|_2(\pi L)^n}{n!} \fstop
\]
These bounds are summable, so the series converges in trace norm;
its kernel converges in $L^2$ to the displayed exponential kernel.
Thus $M_fP$ is trace class, and so is its adjoint $PM_f$.
The same argument applies to every bounded multiplier with this
spatial support.

Having verified the trace-class hypotheses, we now express the moment generating function of the centred statistic through the Fredholm determinant. We then apply the Breuer–Duits estimate and use Cauchy’s formula to obtain the desired cumulant bounds.

\medskip
\noindent\emph{Step 1.}
Define
$$
 \Pi\defi M_{\1_{[-L,L]}}\comma\qquad
 T_z\defi M_{e^{zf}-1}P\comma\qquad
 X_f\defi\sum_{x\in\xi}f(x)\fstop
$$
The preceding argument shows that $T_z$ is trace class.
Since $f$ vanishes outside $[-L,L]$, so does $e^{zf}-1$, 
hence $\Pi T_z=T_z$.
By applying $\det(I+AB)=\det(I+BA)$ with $A=\Pi$ and $B=T_z$, we obtain
$$
 \det_{L^2(\R)}(I+T_z)
 =\det_{L^2(\R)}(I+T_z\Pi)\fstop
$$
The operator $T_z\Pi=\Pi T_z\Pi$ is zero on
$L^2(\R\setminus[-L,L])$ and maps $L^2([-L,L])$ into itself.
Letting
$$
 P_L\defi\left.\Pi P\Pi\right|_{L^2([-L,L])}\comma
$$
its restriction to this latter space is $M_{e^{zf}-1}P_L$.
Therefore
$$
\begin{aligned}
 \det_{L^2(\R)}(I+T_z)
 &=\det_{L^2([-L,L])}\!\bigl(I+M_{e^{zf}-1}P_L\bigr)
 &=\mbbE[e^{zX_f}]\fstop
\end{aligned}
$$
The last equality is Proposition~\ref{prop:fredholm}, applied
on the compact interval containing the support of $f$.

\medskip
\noindent\emph{Step 2.}
The exponential series gives $T_0=0$ and $T'_0=M_fP$
with differentiation in trace norm.
By taking the differentiation for the preceding identity at $z=0$, we have
$$
\begin{aligned}
 \Tr(M_fP)
 &=\left.\frac{\diff}{\diff z}
        \det_{L^2(\R)}(I+T_z)\right|_{z=0} =\left.\frac{\diff}{\diff z}\mbbE[e^{zX_f}]\right|_{z=0}
 =\mbbE[X_f]
 =\int_{\R}f(x)\diff x\fstop
\end{aligned}
$$
Here the last equality uses the unit intensity of the sine point process.
Thus the centring in $\mcN(f)=X_f-\int f$ is  the
linear trace term in the Breuer--Duits estimate.

For $|z|\le1/3$, the operator norm satisfies
$$
 \|T_z\|
 \le e^{|z|\|f\|_\infty}-1
 \le e^{|z|}-1<1\fstop
$$
Hence $I+T_z$ is invertible and its Fredholm determinant is nonzero.
We use the analytic logarithm that vanishes at $z=0$.
Then we obtain that the centred cumulant generating function is 
$$
 \Lambda_f(z)\defi\log\mbbE[e^{z\mcN(f)}]
 =\log\det_{L^2(\R)}(I+T_z)-z\Tr(M_fP)\fstop
$$
By Breuer--Duits~\cite[Lemma~3.2]{BreuerDuits2014}, we have 
\begin{equation}\label{eq:one-time-log-mgf}
 \big|\log\mbbE e^{z\mcN(f)}\big|
 \le C_0|z|^2\|[P,M_f]\|_2^2\comma
 \qquad |z|\le\frac13\fstop
\end{equation}
This estimate holds also for complex $z$.
Indeed, write $z=r e^{i\theta}$ and apply the cited lemma with
the real parameter $r=|z|$ and the bounded operator
$h=e^{i\theta}M_f$.

\medskip
\noindent\emph{Step 3: extract the cumulant bound.}
By~\eqref{eq:cumulant-definition}, the coefficient of $z^m$
in $\Lambda_f(z)$ is $\kappa_m(\mcN(f))/m!$.
The circle of radius $r=1/4$ lies inside its disc of analyticity.
Using Cauchy's coefficient formula, 
$$
 \frac{\kappa_m(\mcN(f))}{m!}
 =\frac1{2\pi i}\int_{|z|=r}
       \frac{\Lambda_f(z)}{z^{m+1}}\diff z\fstop
$$
The contour is traversed counterclockwise and has length $2\pi r$.
 Using~\eqref{eq:one-time-log-mgf},
we obtain
$$
\begin{aligned}
 \frac{|\kappa_m(\mcN(f))|}{m!}
 &\le r^{-m}\max_{|z|=r}|\Lambda_f(z)|
 \le C_0 r^{2-m}\|[P,M_f]\|_2^2
 =C_0\,4^{m-2}\|[P,M_f]\|_2^2\fstop
\end{aligned}
$$
Regarding the choice of constant: for instance, $C=4\max\{1,C_0\}$ satisfies
$C_0\,4^{m-2}\le C^m$ for every $m\ge3$.
By multiplying by $m!$, we obtain~\eqref{eq:sine-cumulant-bound}.
For $f=f_R$, Lemma~\ref{lem:comm-logR} gives
$\|[P,M_{f_R}]\|_2^2\le C_1\log R$.
By replacing $C$ by $C\max\{1,C_1\}$, we complete the proof of  the stated logarithmic bound.
\end{proof}

\subsubsection{Uniform cumulant bounds on finitely many time layers}

\begin{prop}[Uniform finite-time cumulant bound]
\label{prop:cumulant-bound}\label{prop:multitime-cumulant-bound}
There exists a constant $C<\infty$, depending only on the fixed
cut-off $\chi$, with the following property.
For any integer $r\ge1$, real coefficients
$\theta_0,\ldots,\theta_r$, observation times
$0=t_0<t_1<\cdots<t_r$, and $R\ge2$, set
$$
 Y_R\defi\sum_{a=0}^r\theta_a
            \int_\R f_R(x)\xi_{t_a}(\diff x)\comma
 \qquad
\Theta\defi\max_{0\le a\le r}|\theta_a|\fstop
$$
Then, for every integer $m\ge3$,
\begin{equation} \notag %\label{eq:multitime-cumulant-bound}
 |\kappa_m(Y_R)|
 \le m!\bigl(C(r+1) \Theta \bigr)^m\log R\fstop
\end{equation}
In particular, for all $t>0$ and $R\ge2$,
$$
 |\kappa_m(Z_R(t))|\le m!C^m\log R\fstop
$$
Consequently, with $R(t)=t^2$, for every $t\ge2$,
\begin{equation} \notag %\label{eq:cumulant-bound-final}
 |\kappa_m(Z_{R(t)}(t))|\le m!C^m\log t\fstop
\end{equation}
All bounds are uniform in the observation times and their gaps.
\end{prop}

\begin{proof}
Fix an integer $m\ge3$ and put
$$
 f\defi f_R\comma\qquad 
  d_R\defi\|[P,M_f]\|_2\fstop
$$
Since the assertion is immediate if $\Theta=0$, we may assume $\Theta>0$ for the rest of the proof.

\noindent\emph{Step 1: the time-word expansion.}
Use the gauged blocks $K_{t_a,t_b}$ defined in
Section~\ref{subsec:mixed-time-traces} and represented
by~\eqref{eq:multitime-gauge-blocks}. Define
$$
 E\defi\bigsqcup_{a=0}^r(\{t_a\}\times\R)\comma\qquad
 \Xi\defi\bigsqcup_{a=0}^r(\{t_a\}\times\xi_{t_a})\fstop
$$
By~\eqref{eq:extended-sine-correlations} and the gauge
transformation~\eqref{eq:gauge-kernel}, its correlation kernel
is $\widetilde K_{\sin}$. Set $F(t_a,x)\defi\theta_af(x)$.
The linear statistic \(X_F\) appearing in Proposition~\ref{prop:cumulant-trace} is then exactly \(Y_R\):
\[
X_F
=\sum_{(t_a,x)\in\Xi}F(t_a,x)
=\sum_{a=0}^{r}\theta_a\sum_{x\in\xi_{t_a}}f_R(x)
=Y_R.
\]

Set
$$
 \Lambda_R\defi\supp f\comma\qquad
 A\defi\supp_E F
 =\bigsqcup_{\{a:\theta_a\ne0\}}
       (\{t_a\}\times\Lambda_R)\fstop
$$
Recall that the reference measure $\nu\defi\#\times\diff x$ is counting
measure on the finite set of times tensored by
Lebesgue measure in space.
We also recall that  $K_A$ is the integral operator
on $L^2(A,\nu)$ obtained by restricting both arguments of the
gauged kernel to $A$, viz., 
$$
 (K_Ah)(z)\defi
 \int_A\widetilde K_{\sin}(z,w)h(w)\nu(\diff w)\comma
 \qquad z\in A\fstop
$$
Lemma~\ref{lem:local-trace-class} then shows that $K_A$ is a trace class operator.
We may therefore apply the cumulant formula of
Proposition~\ref{prop:cumulant-trace} to $F$.
For a positive composition
$\mathbf p=(p_1,\ldots,p_\ell)$ of $m$, i.e.,  $p_j\ge1$ and $p_1+\cdots+p_\ell=m$, define
$$
 \mathcal T_{\mathbf p}\defi
 \Tr_{L^2(A,\nu)}
 \bigl(M_{F^{p_1}}K_A\cdots M_{F^{p_\ell}}K_A\bigr)
 \fstop
$$
These traces are well defined because the multiplication
operators are bounded and $K_A$ is a trace class opeartor.
Set
$$
 c_{\mathbf p}\defi
 (-1)^{\ell-1}\frac{m!}{\ell\,p_1!\cdots p_\ell!}\fstop
$$
To define the individual time words with the same localisation
convention, put $\Pi_R\defi M_{\1_{\Lambda_R}}$ on $L^2(\R)$.
For an ordered word
$\mathbf a=(a_1,\ldots,a_\ell)\in\{0,\ldots,r\}^{\ell}$,
use the closing convention $a_{\ell+1}=a_1$ and set
$$
 W_{\mathbf a,\mathbf p}\defi
 \Tr_{L^2(\R)}\Big(
 M_{(\theta_{a_1}f)^{p_1}}K_{t_{a_1},t_{a_2}}\cdots
 M_{(\theta_{a_\ell}f)^{p_\ell}}K_{t_{a_\ell},t_{a_1}}
 \Pi_R
 \Big)\fstop
$$

We now apply Lemma~\ref{lem:trace-localisation}
to
$$
 G_j(t_a,x)\defi
 \1_{\{a=a_j\}}(\theta_af(x))^{p_j}\comma
 \qquad 1\le j\le\ell\fstop
$$
Each multiplier selects its prescribed time layer and is
supported in $A$. The lemma identifies the resulting restricted
trace with $W_{\mathbf a,\mathbf p}$ and shows that it is
well defined for every $\ell\ge1$.
For $\ell\ge2$, the final projection
$\Pi_R$ may be omitted by~\eqref{eq:trace-full-space}.
For $\ell=1$, we retain it as in~\eqref{eq:trace-localisation};
every such word is constant.

Put $I\defi\{a\in\{0,\ldots,r\}:\theta_a\ne0\}$ and
$\mcH_\Lambda\defi L^2(\Lambda_R)$, so that
$L^2(A,\nu)\simeq\bigoplus_{a\in I}\mcH_\Lambda$.
Relative to this direct sum, the $(a,b)$-block
$D_{ab}:\mcH_\Lambda\to\mcH_\Lambda$ of a bounded operator
$D$ on $L^2(A,\nu)$ is characterised by
$$
 (Dh)_a=\sum_{b\in I}D_{ab}h_b\comma
 \qquad h=(h_b)_{b\in I}\fstop
$$
Thus column $b$ specifies the input time $t_b$, and row $a$
specifies the output time $t_a$.
Identify $\mcH_\Lambda$ with the subspace of $L^2(\R)$
consisting of functions that vanish outside $\Lambda_R$.
For $a,b\in I$, let $K^R_{ab}$ be the restriction of
$\Pi_RK_{t_a,t_b}\Pi_R$ to $\mcH_\Lambda$.
Writing $D_j\defi M_{F^{p_j}}K_A$, its time-layer blocks are
$$
 (D_j)_{ab}=M_{(\theta_af)^{p_j}}K^R_{ab}\fstop
$$
Therefore, 
$$
 \begin{aligned}
 (D_1\cdots D_\ell)_{a_1a_1}
 =\sum_{a_2,\ldots,a_\ell\in I}
 M_{(\theta_{a_1}f)^{p_1}}K^R_{a_1a_2}
 \cdots
 M_{(\theta_{a_\ell}f)^{p_\ell}}K^R_{a_\ell a_1}\fstop
 \end{aligned}
$$
For $\ell=1$, there are no intermediate indices to sum over.
Taking the trace now sums over the starting layer~$a_1$:
$$
 \begin{aligned}
 \mathcal T_{\mathbf p}
 &=\sum_{a_1\in I}
   \Tr_{\mcH_\Lambda}
   \bigl((D_1\cdots D_\ell)_{a_1a_1}\bigr)
 &=\sum_{\mathbf a\in I^\ell}W_{\mathbf a,\mathbf p}
  =\sum_{\mathbf a\in\{0,\ldots,r\}^{\ell}}
       W_{\mathbf a,\mathbf p}\fstop
 \end{aligned}
$$
Here Lemma~\ref{lem:trace-localisation} identifies each restricted
trace with $W_{\mathbf a,\mathbf p}$, including its final
projection $\Pi_R$. The last equality holds because a word
containing an index with $\theta_a=0$ contributes zero.

Consequently, the cumulant formula~\eqref{eq:cumulant-trace} is
$$
 \begin{aligned}
 \kappa_m(Y_R)
 &=\sum_{\ell=1}^m
\sum_{\substack{p_1+\cdots+p_\ell=m\\p_j\ge1}}
      c_{\mathbf p}\mathcal T_{\mathbf p}
 &=\sum_{\ell=1}^m
\sum_{\substack{p_1+\cdots+p_\ell=m\\p_j\ge1}}
      c_{\mathbf p}
      \sum_{\mathbf a\in\{0,\ldots,r\}^{\ell}}
         W_{\mathbf a,\mathbf p}\fstop
 \end{aligned}
$$
In the above expansion, define
$\kappa_m^{\mrmconst}(Y_R)$ to be the sum of all terms
$c_{\mathbf p}W_{\mathbf a,\mathbf p}$ whose time indices
satisfy $a_1=\cdots=a_\ell$. Define
$\kappa_m^{\mathrm{mix}}(Y_R)$ to be the sum of all remaining
terms, whose words visit at least two distinct time layers.
Each sum includes all lengths and compositions appearing
in the expansion. Thus
$$
\kappa_m(Y_R)
=
\kappa_m^{\mrmconst}(Y_R)
+
\kappa_m^{\mathrm{mix}}(Y_R)\fstop
$$
In the following steps, we evaluate each terms separately.

\medskip
\noindent\emph{Step 2: constant time words.}
If $\mathbf a=(a,\ldots,a)$, every kernel block is $P$.
By factoring out the coefficients, we obtain 
$$
 W_{(a,\ldots,a),\mathbf p}
 =\theta_a^m V_{\mathbf p}\comma
 \qquad
 V_{\mathbf p}\defi
 \Tr_{L^2(\R)}
 \bigl(M_{f^{p_1}}P\cdots M_{f^{p_\ell}}P\Pi_R\bigr)\fstop
$$
As above, the final projection may be omitted when $\ell\ge2$.
By applying Proposition~\ref{prop:cumulant-trace} with the setting of one-time, we obtain 
$$
 \sum_{\ell=1}^m
 \sum_{\substack{p_1+\cdots+p_\ell=m\\p_j\ge1}}
       c_{\mathbf p}V_{\mathbf p}
 =\kappa_m(\mcN(f))\fstop
$$
Here the centring does not change the cumulant because $m\ge3$,
and the stationarity identifies the one-time law on every layer.
Consequently, we have 
$$
 \kappa_m^{\mrmconst}(Y_R)
 =\sum_{a=0}^r\theta_a^m\kappa_m(\mcN(f))\fstop
$$
Using Proposition~\ref{prop:sine-cumulant-bound} and
$\sum_{a=0}^r |\theta_a|^m\le (r+1) \Theta^m$, we obtain
$$
\begin{aligned}
 |\kappa_m^{\mrmconst}(Y_R)|
 &\le m!C_{\mathrm s}^{\,m}d_R^2\sum_{a=0}^r|\theta_a|^m
 &\le m!\, (r+1) (C_{\mathrm s}\Theta)^m d_R^2\fstop
\end{aligned}
$$

\noindent\emph{Step 3: a bound for each mixed word.}
Fix a nonconstant cyclic word $\mathbf a$ and a positive
composition $\mathbf p=(p_1,\ldots,p_\ell)$ of $m$.
Let $T_{\mathbf a,\mathbf p}$ be the operator in
Lemma~\ref{lem:two-interfaces}. By factoring the scalar coefficients
out of the multiplication operators, we have 
$$
 W_{\mathbf a,\mathbf p}
 =\left(\prod_{j=1}^{\ell}\theta_{a_j}^{p_j}\right)
       \Tr(T_{\mathbf a,\mathbf p})\comma\qquad
 \left|\prod_{j=1}^{\ell}\theta_{a_j}^{p_j}\right|
 \le\Theta^{\sum_jp_j}=\Theta^m\fstop
$$
Now Lemma~\ref{lem:two-interfaces} applies. Moreover,
the power estimate~\eqref{eq:comm-power} gives that,  for each $j$, we have 
$$
 \|[P,M_{f^{p_j}}]\|_2
 \le p_j\|[P,M_f]\|_2=p_jd_R\fstop
$$
Consequently,
$$
\begin{aligned}
 |\Tr(T_{\mathbf a,\mathbf p})|
 &\le\left(\max_{1\le j\le\ell}
                   \|[P,M_{f^{p_j}}]\|_2\right)^2 
 &\le\left(\max_{1\le j\le\ell}p_j\right)^2d_R^2
 \le m^2d_R^2\fstop
\end{aligned}
$$
Here the last inequality uses $p_1+\cdots+p_\ell=m$.
By combining this trace bound with the scalar estimate, we have 
$$
 |W_{\mathbf a,\mathbf p}|\le m^2\Theta^m d_R^2\fstop
$$
The uniformity in the observation times follows from
Lemma~\ref{lem:two-interfaces}.

\medskip
\noindent\emph{Step 4: summation and the two-time specialisation.} 
Let $ n =  r+1$. 
At length $\ell$, there are at most $n^\ell$ time words and
$\binom{m-1}{\ell-1}$ positive compositions. The coefficient bound is
$$
 |c_{\mathbf p}|
 =\frac{m!}{\ell\,p_1!\cdots p_\ell!}\le m!\fstop
$$
Since a mixed word has $\ell\ge2$, we obtain
$$
 |\kappa_m^{\mathrm{mix}}(Y_R)|
 \le m!m^2\Theta^m d_R^2
       \sum_{\ell=2}^m n^\ell\binom{m-1}{\ell-1}\fstop
$$
Using the binomial theorem,  
$$
 \sum_{\ell=2}^m n^\ell\binom{m-1}{\ell-1}
 \le n\sum_{q=0}^{m-1}n^q\binom{m-1}{q}
 =n(n+1)^{m-1}\fstop
$$
Using $m^2\le4^m$ and $n(n+1)^{m-1}\le(2n)^m$, we obtain 
$$
 |\kappa_m^{\mathrm{mix}}(Y_R)|
 \le m!(8n\Theta)^m d_R^2\fstop
$$
Adding the constant-time bound and using $n\le n^m$, we find
$$
\begin{aligned}
 |\kappa_m(Y_R)|
 &\le m!\bigl((C_{\mathrm s}n\Theta)^m+(8n\Theta)^m\bigr)d_R^2 
 &\le m!\bigl((C_{\mathrm s}+8)n\Theta\bigr)^m d_R^2\fstop
\end{aligned}
$$
Finally, Lemma~\ref{lem:comm-logR} gives
$$
 d_R^2\le C_\chi\log R\comma
 \qquad
 |\kappa_m(Y_R)|\le m!(C_0n\Theta)^m\log R\comma
$$
where $C_0\defi(C_{\mathrm s}+8)\max\{1,C_\chi\}$.
Both the constants appearing in the one-time and mixed-word bounds are
independent of the  time-points. Hence, the constant $C_0$ depends only on the fixed cut-off function.

For the two-time statistic, we take
$$
 r=1\comma\qquad t_1=t\comma\qquad
 (\theta_0,\theta_1)=(-1,1)\fstop
$$
Then $Y_R=Z_R(t)$, $n=2$ and $\Theta=1$. Hence, 
$$
 |\kappa_m(Z_R(t))|\le m!(2C_0)^m\log R\fstop
$$
With $R(t)=t^2$ and $t\ge2$,
$$
 |\kappa_m(Z_{R(t)}(t))|
 \le2m!(2C_0)^m\log t
 \le m!(4C_0)^m\log t\fstop
$$
By taking $C=4C_0$, we obtain all the stated bounds with the same constant.
\end{proof}

\subsection{Proof of the CLT}

\begin{proof}[Proof of Theorem~\ref{thm:CLT}]
We define $X_t\defi Z_{R(t)}(t)$ with $R(t)=t^2$.
Recall that $Z_R(t)=\mcN_t(f_R)-\mcN_0(f_R)$ with 
$f_R(x)=\chi(x/R)\1_{\{x>0\}}$.
By Lemma~\ref{lem:CLT-reduction}, we have 
$$
 \delta_t\defi\|X_t-J(t)\|_{L^2}=O(t^{-1/2})\fstop
$$
We first compare their variances, then prove the CLT for $X_t$.
We then transfer it to $J(t)$.

\medskip
\noindent\emph{Step 1: comparison of the variances.}
Recalling that both variables are centred, we have that their standard deviations are equal
 to their $L^2$ norms. The reverse triangle inequality then gives
$$
 \big|\sqrt{\Var(X_t)}-\sqrt{\Var(J(t))}\big|
 =\big|\|X_t\|_{L^2}-\|J(t)\|_{L^2}\big|
 \le\delta_t\fstop
$$
In particular, $\sqrt{\Var(X_t)}\le\sqrt{\Var(J(t))}+\delta_t$.
By the factorisation of  the difference of the variances, we obtain
$$
\begin{aligned}
 |\Var(X_t)-\Var(J(t))|
 &=\big|\sqrt{\Var(X_t)}-\sqrt{\Var(J(t))}\big|
   \bigl(\sqrt{\Var(X_t)}+\sqrt{\Var(J(t))}\bigr)\\
 &\le\delta_t\bigl(2\sqrt{\Var(J(t))}+\delta_t\bigr)
 =O\!\left(\sqrt{\frac{\log t}{t}}\right)\longrightarrow0\fstop
\end{aligned}
$$
Here the last estimate uses Theorem~\ref{thm:main}.
Consequently, the same theorem also implies
\begin{equation}\label{eq:Var-Xt}
 \sigma_t^2\defi\Var(X_t)
 =\frac{1}{\pi^2}\log t+O(1)\comma
 \qquad\sigma_t^2\asymp\log t\fstop
\end{equation}

\medskip
\noindent\emph{Step 2: the CLT for the compact statistic.}
For sufficiently large $t$, set $Y_t\defi X_t/\sigma_t$.
Then $\kappa_1(Y_t)=0$ as $X_t$ is centred, and $\kappa_2(Y_t)=1$ as $\kappa_2$ is the variance by definition.
For each fixed $m\ge3$, Proposition~\ref{prop:cumulant-bound}
and~\eqref{eq:Var-Xt} give
$$
 |\kappa_m(Y_t)|
 =\frac{|\kappa_m(X_t)|}{\sigma_t^m}
 \le m!C^m\frac{\log t}{\sigma_t^m}
 =O_m\bigl((\log t)^{1-m/2}\bigr)\longrightarrow0\fstop
$$
Consequently, all moments converge to the standard Gaussian
moments. The method of moments therefore gives
$Y_t\Rightarrow\mcN(0,1)$; see
\cite[Example~30.1 and Theorem~30.2]{Billingsley1995}.

\medskip
\noindent\emph{Step 3: transfer to the current.}
The variance comparison gives
$\sigma_t/\sqrt{\Var(J(t))}\to1$, while
Lemma~\ref{lem:CLT-reduction} gives
$$
 \left\|\frac{J(t)-X_t}{\sqrt{\Var(J(t))}}\right\|_{L^2}
 =\frac{\delta_t}{\sqrt{\Var(J(t))}}\longrightarrow0\fstop
$$
Using the decomposition
$$
 \frac{J(t)}{\sqrt{\Var(J(t))}}
 =\frac{\sigma_t}{\sqrt{\Var(J(t))}}Y_t
  +\frac{J(t)-X_t}{\sqrt{\Var(J(t))}}\comma
$$
Slutsky's theorem proves the asserted CLT.
\end{proof}

\section{CLT for the tagged particle}\label{subsec:tracer}

We deduce a central limit theorem for the displacement of a tagged
particle from Theorem~\ref{thm:CLT}. 

\medskip
\noindent\textbf{Labelling convention.}
Recall the labelling convention fixed in the introduction:
\begin{equation}\label{eq:label-convention}
I_0(0)=0\comma \qquad\text{equivalently}\qquad
X_0(0)\le 0 < X_1(0)
\ \ \text{a.s.}
\end{equation}
Recall also the displacement
$$
\Delta(t)\defi X_0(t)-X_0(0)\fstop
$$

\begin{cor}[CLT for the tagged particle]\label{cor:tracer-CLT}
As $t\to\infty$,
$$
\frac{\pi\,(X_0(t)-X_0(0))}{\sqrt{\log t}}
\ \Longrightarrow\ \mcN(0,1)\fstop
$$
\end{cor}

\subsection{Signed counting function and an exact identity}

For a locally finite configuration $\xi=\sum_i\delta_{x_i}$ on $\R$ and $x\in\R$, define the signed counting function
\begin{equation} \notag %\label{eq:signed-count}
N_\xi(x)\defi 
\begin{cases}
\xi((0,x]) ,& x\ge 0,\\[1mm]
-\xi((x,0]) ,& x<0.
\end{cases}
\end{equation}
For  $\xi_t\defi \sum_i\delta_{X_i(t)}$, we write $N_t(x)\defi N_{\xi_t}(x)$.

\begin{lem}\label{lem:N-index}
For every $t\ge 0$ and $x\in\R$,
\begin{equation}\label{eq:N-index}
N_t(x)=I_x(t)-I_0(t)\fstop
\end{equation}
In particular, under the labelling convention \eqref{eq:label-convention},
\begin{equation}\label{eq:J-equals-N}
J(t)=N_t(X_0(t))\fstop
\end{equation}
\end{lem}

\begin{proof}
If $x\ge 0$, then $I_x(t)-I_0(t)$ counts the number of indices $i$ such that $0<X_i(t)\le x$, i.e., \ $\xi_t((0,x])$.
If $x<0$, then $I_x(t)-I_0(t)$ is the negative of the number of indices with $x<X_i(t)\le 0$, i.e., \ $-\xi_t((x,0])$.
This proves \eqref{eq:N-index}.

At $x=X_0(t)$, we have $I_{X_0(t)}(t)=0$.
As the labelling convention gives $J(t)=-I_0(t)$, 
\eqref{eq:N-index} becomes~\eqref{eq:J-equals-N}.
\end{proof}

\subsection{\texorpdfstring{A uniform counting-error bound for sine$_2$}{A uniform counting-error bound for sine2}}

\begin{lem}\label{lem:number-variance}
Let $\xi\sim\sine_2$ have intensity $1$. There exists $C<\infty$ such that, for all $L\ge 2$,
$$
\Var\big(\xi([0,L])\big)\le C\log L\fstop
$$
The same bound holds for $\xi_t$ for every $t$ under the stationary law.
\end{lem}

\begin{proof}
Due to the asymptotic of the number varianec of the sine point process, we have 
$$
 \Var(\xi([0,L]))\sim\frac1{\pi^2}\log L
 \qquad (L\to\infty)\comma
$$
by \cite[Eq.~(3.15)]{Soshnikov2000DRPF}, taking the single
Fourier interval $[-\pi,\pi]$ and translating the counting
interval. This gives the bound for large $L$. For the remaining
bounded range $L\ge2$, use
$\Var(\xi([0,L]))\le\mbbE[\xi([0,L])]=L$ and take a sufficiently large $C$.
The stationarity of $\xi_t$ gives the same conclusion for $\xi_t$ at every
fixed $t$.
\end{proof}

\subsection{Tagged particle displacement}
We write $\xrightarrow{\mbbP}$ for convergence in probability. 
\begin{prop}\label{prop:tracer-minus-current}
As $t\to\infty$,
$$
\frac{X_0(t)-X_0(0)-J(t)}{\sqrt{\log t}}\xrightarrow{\mbbP}0\fstop
$$
\end{prop}

%\begin{proof}
\begin{proof}
Due to the identity $J(t)=N_t(X_0(t))$, we have 
\begin{equation} \notag %\label{eq:tracer-current-diff}
 X_0(t)-J(t)=X_0(t)-N_t(X_0(t))\fstop
\end{equation}
We therefore estimate the counting discrepancy uniformly in terms of radius of   interval. We then show that the tagged
particle belongs to this interval with probability tending to one.

Set
\begin{align}\label{:44a}&
 M(t)\defi(\log t)^{3/4}\comma
 \qquad u(t)\defi(\log t)^{9/20}\fstop
\end{align}
It holds that 
\begin{align}\label{:44b}&
 \frac{\sqrt{\log t}}{M(t)}\longrightarrow0\comma
 \qquad
 \frac{u(t)}{\sqrt{\log t}}\longrightarrow0\comma
 \qquad
 \frac{M(t)\log M(t)}{u(t)^2}\longrightarrow0\fstop
\end{align}

\medskip
\noindent\emph{Step 1.}
For every nonzero integer $n$, by the intensity-one property and
the spatial translation invariance, we have 
\begin{align}\label{:44c}&
 \mbbE[N_t(n)]=n\comma
 \qquad
 \Var(N_t(n))\le C\log(2+|n|)\fstop
\end{align}
For $n<0$, this follows by applying the number-variance estimate
to $N_t(n)=-\xi_t((n,0])$.
Lemma~\ref{lem:number-variance} applies when $|n|\ge2$ while
the bound for $|n|=1$ follows from
$\Var(\xi_t(I))\le\mbbE[\xi_t(I)]$.
%Deterministic endpoints do not affect the counts almost surely,
%since the intensity measure has no atoms. Also, $N_t(0)=0$.

Chebyshev's inequality and \eqref{:44c} yield 
$$
 \mbbP\bigl(|N_t(n)-n|\ge u(t)\bigr)
 \le\frac{C\log(2+|n|)}{u(t)^2}\fstop
$$
There are integers of order $O(M(t))$  such that 
$|n|\le\lceil M(t)\rceil+1$. By summing these bounds and using \eqref{:44b}, we have  
\begin{equation} \label{eq:grid-union}
 \mbbP\left(
  \max_{\substack{n\in\Z\\|n|\le\lceil M(t)\rceil+1}}
       |N_t(n)-n|\ge u(t)
 \right)
 \le\frac{CM(t)\log M(t)}{u(t)^2}\longrightarrow0\fstop 
\end{equation}

To extend this estimate from integers to real numbers, use the monotonicity
of $N_t$. For $n\le x\le n+1$,
$$
 N_t(n)\le N_t(x)\le N_t(n+1)\comma
$$
and hence
$$
 N_t(n)-n-1
 \le N_t(x)-x
 \le N_t(n+1)-(n+1)+1\fstop
$$
Every point of $[-M(t),M(t)]$ lies between two integers in the
preceding grid. Consequently,
\begin{equation} \label{eq:sup-control}
 \sup_{|x|\le M(t)}|N_t(x)-x|
 \le1+\max_{\substack{n\in\Z\\|n|\le\lceil M(t)\rceil+1}}
              |N_t(n)-n|\comma
\end{equation}

By \eqref{eq:grid-union} and \eqref{eq:sup-control}, we obtain 
\begin{equation} \notag %\label{eq:sup-prob}
 \mbbP\left(
  \sup_{|x|\le M(t)}|N_t(x)-x|\ge1+u(t)
 \right)\longrightarrow0\fstop
\end{equation}

\medskip
\noindent\emph{Step 2.}
Define
$$
 \mcE_t\defi
 \left\{\sup_{|x|\le M(t)}|N_t(x)-x|\le2+2u(t)\right\}
 \cap\left\{|J(t)|\le M(t)/2\right\}\fstop
$$
The first event has probability tending to one by Step~1.
Due to Theorem~\ref{thm:main} and Chebyshev's inequality, we have
$$
 \mbbP\bigl(|J(t)|>M(t)/2\bigr)
 \le\frac{4\Var(J(t))}{M(t)^2}
 =O\bigl((\log t)^{-1/2}\bigr)\longrightarrow0\fstop
$$
Thus $\mbbP(\mcE_t)\to1$.

On $\mcE_t$, the endpoint counts satisfy
$$
\begin{aligned}
 N_t(M(t))&\ge M(t)-2-2u(t)>M(t)/2\comma\\
 N_t(-M(t))&\le-M(t)+2+2u(t)<-M(t)/2
\end{aligned}
$$
for all sufficiently large $t$ because of  $u(t)=o(M(t))$.

If $X_0(t)>M(t)$, by the monotonicity and the current identity, we would have 
$$
 J(t)=N_t(X_0(t))\ge N_t(M(t))>M(t)/2\fstop
$$
Similarly, $X_0(t)<-M(t)$ would imply
$$
 J(t)=N_t(X_0(t))\le N_t(-M(t))<-M(t)/2\fstop
$$
Either conclusion contradicts the bound \(|J(t)|\le M(t)/2\) imposed in the definition of \(\mathcal E_t\).
%Both conclusions contradict the definition of $\mcE_t$.
Therefore
$$
 |X_0(t)|\le M(t)\qquad\text{on }\mcE_t\fstop
$$

\medskip
\noindent\emph{Step 3.}
The uniform estimate from Step~1 can now be evaluated at $X_0(t)$.
On $\mcE_t$, we have 
$$
 |X_0(t)-J(t)|
 =|X_0(t)-N_t(X_0(t))|
 \le2+2u(t)\fstop
$$
Since $(2+2u(t))/\sqrt{\log t}\to0$ by \eqref{:44b}, for every $\varepsilon>0$
and all sufficiently large $t$,
$$
 \mbbP\left(
  \frac{|X_0(t)-J(t)|}{\sqrt{\log t}}>\varepsilon
 \right)
 \le\mbbP(\mcE_t^c)\longrightarrow0\fstop
$$
It follows that
\begin{equation}\notag %\label{eq:XminusJ-op}
 \frac{X_0(t)-J(t)}{\sqrt{\log t}}
 \xrightarrow{\mbbP}0\fstop
\end{equation}

Thus,  
$$
 \frac{X_0(t)-X_0(0)-J(t)}{\sqrt{\log t}}
 \xrightarrow{\mbbP}0\comma
$$
as required.
\end{proof}

\begin{proof}[Proof of Corollary~\ref{cor:tracer-CLT}]
By Proposition~\ref{prop:tracer-minus-current}, Theorem~\ref{thm:CLT} and Slutsky's theorem, 
 the asserted is proved. 
\end{proof}

\section{Finite-dimensional process CLT}\label{sec:process-clt}

We now extend the one-time CLT to the multi-scaling limit.  
Throughout this section, $\alpha\downarrow0$ is the time-scaling parameter,
and $T_\alpha\defi\log(1/\alpha)$. We use $Z_R$ as in Section~\ref{sec:CLT}.

\begin{lem}\label{lem:current-covariance}
Let $V(t)\defi \Var(J(t))$.  For $s,t>0$ with $s\ne t$,
\begin{equation}\label{eq:current-covariance-identity}
\Cov(J(s),J(t))
=\frac12\{V(s)+V(t)-V(|t-s|)\}\comma
\end{equation}
and $\Cov(J(t),J(t))=V(t)$.  Consequently, for fixed $u,v>0$,
\begin{equation}\notag %\label{eq:hom-current-cov-limit}
\frac{\pi^2}{T_\alpha}\,
\Cov\big(J(u/\alpha),J(v/\alpha)\big)
\longrightarrow
\begin{cases}
1, & u=v,\\[1mm]
\dfrac12, & u\ne v.
\end{cases}
\end{equation}
\end{lem}

\begin{proof}
For $0<s<t$,
$J(t)-J(s)=I_0(s)-I_0(t)$ is  the current through the origin during $[s,t]$.
By the stationarity, it has the same law as $J(t-s)$, hence
$\Var(J(t)-J(s))=V(t-s)$.  
Formula~\eqref{eq:current-covariance-identity} follows just by expanding   the variance formula, and  the case $s>t$ is the same as it is symmetric.
For fixed $c>0$, due to Theorem~\ref{thm:main}, we have 
$V(c/\alpha)=\pi^{-2}T_\alpha+O_c(1)$.
On the diagonal, this proves the limit $1$, and 
off the diagonal, by substituting it for the three variances
in~\eqref{eq:current-covariance-identity}, we obtain that  the limit is  $1/2$.
\end{proof}

\begin{lem}\label{lem:simultaneous-compact-approx}
Fix $u_1,\ldots,u_n>0$ and
$\theta_1,\ldots,\theta_n\in\R$.  With $R_\alpha\defi \alpha^{-2}$,
\begin{equation}\label{eq:simultaneous-compact-approx}
\frac{1}{\sqrt{T_\alpha}}
\sum_{j=1}^n\theta_j\big(J(u_j/\alpha)-Z_{R_\alpha}(u_j/\alpha)\big)
\xrightarrow[L^2]{\alpha\downarrow0}0\fstop
\end{equation}
\end{lem}

\begin{proof}
By applying  Proposition~\ref{prop:direct-current-L2} at
$t=u_j/\alpha$ with $R=R_\alpha=\alpha^{-2}$, we have 
$$
\|J(u_j/\alpha)-Z_{R_\alpha}(u_j/\alpha)\|_{L^2}
\le\|\psi\|_{L^2}\sqrt{u_j\alpha}\fstop
$$
The triangle inequality therefore shows
$$
\left\|\sum_{j=1}^n\theta_j
\bigl(J(u_j/\alpha)-Z_{R_\alpha}(u_j/\alpha)\bigr)\right\|_{L^2}
\le\|\psi\|_{L^2}\sqrt\alpha
\sum_{j=1}^n|\theta_j|\sqrt{u_j}\longrightarrow0\fstop
$$
Dividing by $\sqrt{T_\alpha}$,  
\eqref{eq:simultaneous-compact-approx} is shown.
\end{proof}

\begin{thm}\label{thm:current-process-clt}
For $u>0$, define
$$
\mcJ_\alpha(u)\defi \frac{\pi J(u/\alpha)}{\sqrt{\log(1/\alpha)}}\fstop
$$
Then, for every finite collection $u_1,\ldots,u_n>0$,
$$
(\mcJ_\alpha(u_1),\ldots,\mcJ_\alpha(u_n))
\Longrightarrow
(G(u_1),\ldots,G(u_n))\comma \qquad \alpha\downarrow0\comma
$$
where $G$ is the centred Gaussian field with covariance
\begin{equation}\label{eq:homogenized-current-covariance}
\mbbE[G(u)G(v)]
=\frac12+\frac12\1_{\{u=v\}}
=\begin{cases}
1,&u=v,\\[1mm]
\dfrac12,&u\ne v.
\end{cases}
\end{equation}
\end{thm}

\begin{proof}
By the Cram\'er--Wold theorem, it suffices to consider
$L_\alpha\defi \sum_{j=1}^n\theta_j\mcJ_\alpha(u_j)$. 
Without loss of generality, we may assume
$0<u_1<\cdots<u_n$.  By Lemma~\ref{lem:simultaneous-compact-approx}, we may replace
$L_\alpha$  by
$$
\widetilde L_\alpha\defi 
\frac{\pi}{\sqrt{T_\alpha}}
\sum_{j=1}^n\theta_j Z_{R_\alpha}(u_j/\alpha)\comma \qquad R_\alpha=\alpha^{-2}\fstop
$$
The numerator is a  linear statistic at times
$0,u_1/\alpha,\ldots,u_n/\alpha$ with coefficients
$\theta_0=-\sum_{j=1}^n\theta_j$ at time $0$ and $\theta_j$ at time
$u_j/\alpha$.  Due to Proposition~\ref{prop:multitime-cumulant-bound},  for every
fixed $m\ge3$, we have 
$$
|\kappa_m(\widetilde L_\alpha)|
\le \frac{\pi^m}{T_\alpha^{m/2}}C_{m,n,\theta}\log R_\alpha
=O(T_\alpha^{1-m/2})\to0\fstop
$$
The first cumulant is zero. By Lemma~\ref{lem:current-covariance},
the variance of $L_\alpha$ converges to
$$
 v_\theta\defi
 \frac12\left(\sum_{j=1}^n\theta_j\right)^2
 +\frac12\sum_{j=1}^n\theta_j^2\comma
$$
the quadratic form associated with~\eqref{eq:homogenized-current-covariance}
at the distinct times $u_1,\ldots,u_n$. Both statistics are centred, so
$$
 \left|\sqrt{\Var(\widetilde L_\alpha)}
       -\sqrt{\Var(L_\alpha)}\right|
 \le\|\widetilde L_\alpha-L_\alpha\|_{L^2}\longrightarrow0\fstop
$$
Thus the second cumulant of $\widetilde L_\alpha$ also tends to
$v_\theta$. The moment--cumulant formula, as in the proof of
Theorem~\ref{thm:CLT}, gives
$\widetilde L_\alpha\Longrightarrow\mcN(0,v_\theta)$.
If all coefficients vanish, this conclusion is immediate.
Slutsky's theorem and the Cram\'er--Wold theorem complete the proof.
\end{proof}

\begin{thm}[Tagged-particle process: multi-time scaling limit]\label{thm:tagged-process-clt}
Let
$$
\Delta(t)\defi X_0(t)-X_0(0)\comma \qquad
\mcX_\alpha(u)\defi \frac{\pi\Delta(u/\alpha)}{\sqrt{\log(1/\alpha)}}\fstop
$$
Then, for every finite collection $u_1,\ldots,u_n>0$,
$$
(\mcX_\alpha(u_1),\ldots,\mcX_\alpha(u_n))
\Longrightarrow
(G(u_1),\ldots,G(u_n))\comma \qquad \alpha\downarrow0\comma
$$
where $G$ has covariance \eqref{eq:homogenized-current-covariance}.
\end{thm}

\begin{proof}
For each fixed $u>0$, applying Proposition~\ref{prop:tracer-minus-current} with
$t=u/\alpha$, we have 
$$
\frac{\Delta(u/\alpha)-J(u/\alpha)}{\sqrt{\log(1/\alpha)}}
\xrightarrow{\mbbP}0\comma
$$
because $\log(u/\alpha)\sim\log(1/\alpha)$. 
For fixed $u_1,\ldots,u_n>0$ and every $\varepsilon>0$,
$$
\begin{aligned}
&\mbbP\left(
 \max_{1\le j\le n}
 \frac{|\Delta(u_j/\alpha)-J(u_j/\alpha)|}
      {\sqrt{\log(1/\alpha)}}>\varepsilon
\right)
\le
 \sum_{j=1}^n
 \mbbP\left(
 \frac{|\Delta(u_j/\alpha)-J(u_j/\alpha)|}
      {\sqrt{\log(1/\alpha)}}>\varepsilon
 \right)
 \xrightarrow{\alpha \downarrow 0} 0\fstop
\end{aligned}
$$
 The desired result follows from Theorem~\ref{thm:current-process-clt}
and Slutsky's theorem.
\end{proof}

\begin{rem}[No usual path-space topology]\label{rem:no-path-topology}
The convergence in Theorems~\ref{thm:current-process-clt} and
\ref{thm:tagged-process-clt} is finite-dimensional only.  The limiting field does not have either  continuous or cadlag paths: for every $u\ne v$,
$$
\mbbE[(G(u)-G(v))^2]=1\fstop
$$
Thus the theorem should not be read as convergence in $C([a,b])$ or $D([a,b])$
for $0<a<b<\infty$.
\end{rem}

\subsection*{Conflict of Interest}
The authors have no conflict of interest for this research. 

\subsection*{Data availability statement}
No data was generated during this research.

\end{document}